\documentclass[12pt]{amsart}
\usepackage{latexsym}
\usepackage{amssymb}
\usepackage{amsfonts}
\usepackage{epsfig}
\usepackage{amsmath}
\usepackage{amscd}
\usepackage{psfrag}
\usepackage[all]{xy}

\newtheorem{defn0}{Definition}[section]
\newtheorem{prop0}[defn0]{Proposition}
\newtheorem{thm0}[defn0]{Theorem}
\newtheorem{lemma0}[defn0]{Lemma}
\newtheorem{corollary0}[defn0]{Corollary}
\newtheorem{example0}[defn0]{Example}
\newtheorem{conjecture0}[defn0]{Conjecture}
\newtheorem{notation0}[defn0]{Notation}

\theoremstyle{remark}
\newtheorem{remark0}[defn0]{Remark}

\newenvironment{defn}{\begin{defn0}}{\end{defn0}}
\newenvironment{prop}{\begin{prop0}}{\end{prop0}}
\newenvironment{thm}{\begin{thm0}}{\end{thm0}}
\newenvironment{lem}{\begin{lemma0}}{\end{lemma0}}
\newenvironment{cor}{\begin{corollary0}}{\end{corollary0}}
\newenvironment{example}{\begin{example0}\rm}{\end{example0}}
\newenvironment{rem}{\begin{remark0}\rm}{\end{remark0}}

\newenvironment{notation}{\begin{notation0}\rm}{\end{notation0}}

\renewenvironment{proof}{\noindent {\textsc{Proof.}}}{$\square$ \vspace{3mm}}

\newcommand{\D}{\mathbf{\Delta}}
\newcommand{\C}{{\mathbb C}}

\newcommand{\Q}{{\mathbf Q}}

\newcommand{\OO}{{\mathcal O}}

\newcommand{\ZZ}{\mathbb{Z}}

\newcommand{\RR}{\mathbb{R}}

\newcommand{\CC}{\mathbb{C}}

\newcommand{\NN}{\mathcal{N}}

\newcommand{\lra}{\longrightarrow}

\newcommand{\srel}{\stackrel}
\renewcommand{\deg}{\operatorname{deg}}

\newcommand{\eps}{\varepsilon}

\newif\ifprivate
\privatetrue

\def\???{\ifprivate {\bf {???}} \marginpar{{\Huge {\bf ?}}}
\else \fi}

 \numberwithin{equation}{section}

\begin{document}

\title{Geometry of the Kimura 3-parameter model}

\author{Marta Casanellas} \thanks{First author is partially supported by Ministerio de Educaci\'on y Ciencia, programa Ram\'on y Cajal and MTM2006-E14234-C02-02.}
\author{Jes\'{u}s Fern\'{a}ndez-S\'{a}nchez} \thanks{Second author is partially supported by Ministerio de Educaci\'on y Ciencia, programa Juan de la Cierva and MTM2006-E14234-C02-02, MTM2005-01518.}


\maketitle

\begin{abstract}
The Kimura 3-parameter model on a tree of $n$ leaves is one of the
most used in phylogenetics. The affine algebraic variety $W$
associated to it is a toric variety. We study its geometry and we
prove that it is isomorphic to a geometric quotient of the affine
space by a finite group acting on it. As a consequence, we are
able to study the singularities of $W$ and prove that the
biologically meaningful points are smooth points. Then we give an
algorithm for constructing a set of minimal generators of the
localized ideal at these points, for an arbitrary number of leaves
$n$. This leads to a major improvement of phylogenetic reconstruction
methods based on algebraic geometry.
\end{abstract}

\section{Introduction}
The goal of phylogenetic algebraic
geometry is to translate the knowledge
of algebraic geometry into new tools
for phylogenetic inference problems.
The dictionary used in this translation
is based on algebraic statistics, which
allows viewing statistical evolutionary
models as algebraic varieties. The
first approaches in this direction are
due to Allman and Rhodes
\cite{Allman2003} and Pachter and
Sturmfels \cite{Pachter2004}. Since
then, many other authors have
contributed to the development of
phylogenetic algebraic geometry, either
from the more geometric point of view
(see for instance \cite{PhyloAG},
\cite{Sturmfels2005},
\cite{Allman2004b}, \cite{Wisniewski},
\cite{CS}) or from the
applied standing point
(\cite{Evans1993}, \cite{Eriksson05},
\cite{CFS}). The base of algebraic statistics for computational biology were finally set up in the book \cite{ASCB2005}.

The applications of algebraic geometry to phylogenetics rely on
the computation  of the generators of the ideal of the algebraic
variety associated to a statistical evolutionary model on a
phylogenetic tree $T$. In phylogenetics, these generators are
called \textit{phylogenetic invariants} as they are useful to infer
the topology of the tree $T$ (note that in phylogenetics, topology
refers to the topology of the graph $T$ with labels at the
leaves). Phylogenetic invariants have been given for some
algebraic evolutionary models, namely the general Markov model
(\cite{Allman2004b}) and group-based models (Kimura models
\cite{Kimura1980}, \cite{Kimura1981}, and Jukes-Cantor model
\cite{JC69}).

In this paper, we deal with the Kimura 3-parameter model.
As it was shown in \cite{Sturmfels2005}, the Kimura 3-parameter
model on a tree of $n$ species is a toric variety in a
suitable coordinate system (Fourier coordinates). Sturmfels and
Sullivant gave an algorithm to construct a set of generators of
the ideal of this variety for any number of species $n$. For
example, for four species a set of minimal generators contains
8002 binomials of degrees 2, 3 and 4. In a previous paper, we proved that
this set of binomials can be successfully used for phylogenetic
inference (see \cite{CFS}). However, this is a large number of
generators if one considers that the codimension of the variety is
48. Moreover, as  the number of species increases, the codimension
increases exponentially but the number of generators given in
\cite{Sturmfels2005} increases more than exponentially. This makes phylogenetic
reconstruction methods based on this set of generators unfeasible
for larger trees.

The main goal here is to prove that the points of biological interest
are smooth points of the algebraic variety and to provide the
generators of a local complete intersection at these points. To
this end, we prove that this \textit{Kimura variety} $W$ is
isomorphic to the quotient of a certain affine space under the
action of a finite group (Corollary \ref{finit_group}). This
result allows the study of the singular locus of $W$ and shows that
there are no singularities in the points with biological meaning
(Corollary \ref{singular_points}).
We use this result in section 4, where we provide a recursive procedure to give a minimal
sequence of generators for the variety $W$ near these points
(Theorem \ref{invariantsCI}). As an example, the
whole list of these generators in the case of trees with 4 leaves
is given (Example \ref{4leavestree}).

In the paper \cite{Steel1993} the authors also provided a local complete intersection for the Kimura 3-parameter model (they called it a \textit{complete collection of invariants}). In their case the degree of the generators increases exponentially on the number of leaves $n$ and this makes it unfeasible to be used in a phylogenetic reconstruction method for large trees. Our set of generators for the local complete intersection consists of binomials of degree 2 and 4 for any number of leaves $n$ and leads to some hope for the generalization of the method given in \cite{CGS} to arbitrary trees. It is worth mentioning that Hagedorn \cite{Hagedorn2000} also realized that there exists an open set of the variety in which it is sufficient to consider a local complete intersection (this is clear if one knows that the set of singular points on a variety form a Zarisky closed subset). However he did not specify  the open subset nor the set of generators.

This paper is organized as follows. In section 2 we review the
relation between algebraic geometry and statistical evolutionary
models for phylogenetic inference. In this section as well, we
recall the discrete Fourier transform (or Hadamard conjugation)
introduced by Evans and Speed (see \cite{Evans1993}) as a linear
change of coordinates which diagonalizes group-based models. Then
we introduce the algebraic varieties we are interested in and we
set up notation used in the sequel. Section 3 is devoted to the
global study of the geometry of the Kimura variety and to determine its
singular points. In section 4 we perform a local study of the
variety at the biological meaningful points and we give an
algorithm to obtain the generators of a local complete
intersection at these points.

\textbf{Acknowledgments:} The first
author would like to thank L. Pachter
and B. Sturmfels for introducing her to
this subject and encouraging her to
work on it. The second author is deeply
grateful to M. Casanellas for her warm
wellcome on this topic and for giving
him the oportunity of working together.

\section{Preliminaries}

Let $T$ be a tree (i.e. a connected undirected acyclic graph) of
$n$ leaves labelled as $1,2, \dots, n$.
The degree of a node in $T$ is defined as the number of edges
incident to it. Nodes of degree one are leaves $L(T)$, while the
others are internal nodes $N(T)$.
We assume that our trees are trivalent, i.e. internal nodes have
degree 3,
and we call $E(T)$ the set of edges in $T$.
An edge in $T$ is said to be terminal if it contains one leaf. We
write $e_l$ for the terminal edge ending at leaf $l$, $l \in \{1,
\dots, n \}$. It is easy to see that
the number of internal nodes is
$|N(T)|=n-2$ and the number of edges is
$|E(T)|=2n-3$.

\subsection{Algebraic evolutionary models}\label{subsecbio}
In phylogenetics, a tree  $T$ represents the ancestral
relationships (edges) among a set of species (nodes). The leaves
$L(T)$ represent the current species whose phylogenetic history we
wish to infer. The input data is an alignment of $n$ sequences in
the alphabet
 $\Sigma:=\{A,C,G,T\}$ (representing nucleotides)
of length $N$, and one needs to infer the correct phylogenetic
tree that produced the observed alignment.

In order to explain the relationship between phylogenetic
inference and algebraic geometry it is useful to assume for the
moment that the tree is rooted. That is, the graph $T$ is directed
and it has a unique not trivalent node called the \emph{root} of the tree with two edges
emerging from it. This assumption will be removed in subsection
\ref{Kimuravar}.

From the biological standing point, Kimura 3-parameter model is a
stationary Markov model of evolution.  Kimura \cite{Kimura1981}
proposed a statistical model of evolution under the following
assumptions: all sites in the $n$ sequences evolve equally and
independently, the distribution of nucleotides at the root is
uniform and the tree is stationary (and hence all nodes of the
tree have uniform distribution of nucleotides), the evolution of a
species depends only of the node immediately preceding it,
mutations occur randomly and with strictly positive probabilities,
and transitions (mutations between purines A,T or between
pirimydines C,G) occur more often than transversions (mutations
between purines and pirimydines). As all sites evolve
independently and in the same way, one restricts the model to one
site. We describe here an algebraic version of this model (see the
books \cite{ASCB2005} and \cite{ARbook} for an introduction to the
algebraic versions of evolutionary models).

In statistical evolutionary models, to
each node $v$ of the tree $T$ we
associate a discrete random variable
$X_v$ that takes values on
$\Sigma=\{A,C,G,T\}$. The parameters of
algebraic evolutionary models are the
substitution probabilities of
nucleotides along each edge. These
parameters are written in a matrix
indexed by the alphabet elements
$\Sigma$ so that the matrix $S^e$
associated to the edge $e$ is

$ \hspace*{43mm} { A} \hspace*{19mm} { C} \hspace*{18mm} { G}
\hspace*{18mm}  { T}$
$$S^e=
\begin{array}{c}
{ A} \\
{ C} \\
{ G} \\
{ T}
\end{array}
\left(
\begin{array}{cccc}
P( A| A,e) & P( C| A,e) & P({ G}|{ A},e) & P({ T}|{ A},e) \\
P({ A}|{ C},e) & P({ C}|{ C},e) & P({ G}|{ C},e) &P({ T}|{ C},e) \\
P({ A}|{ G},e) & P({ C}|{ G},e) & P({ G}|{ G},e) & P({ T}|{ G},e) \\
P({ A}|{ T},e) & P({ C}|{ T},e) & P({ G}|{ T},e) & P({ T}|{ T},e)
\end{array}
\right)
$$
where $S^e_{x,y}=P(x \mid y , e)$ is
the probability that nucleotide $y$ at
the parent of edge $e$, $s(e)$, mutates
to nucleotide $x$ at the descendant
node $t(e)$. Then the  probability of
observing nucleotides $x_1 \dots x_n$
at the leaves of the tree is given by a
Markov process:

\begin{equation}\label{prob}
p_{x_1 \dots x_n}= \sum_{\{(x_v)_{v \in N(T)} \mid x_v \in \Sigma
\}} \prod_{e \in E(T)} S^e_{x_{s(e)},x_{t(e)}}
\end{equation}
where we assume that if $e=e_l$ is a
terminal edge, then $x_{t(e)}=x_l$.
In
the Kimura 3-parameter model the
substitution matrices have the
following form
$$S^e=\left(
\begin{array}{cccc}
a^e & b^e & c^e & d^e \\
b^e & a^e & d^e & c^e \\
c^e & d^e & a^e & b^e \\
d^e & c^e & b^e & a^e \\
\end{array}
\right)
$$
where $a^e+b^e+c^e+d^e=1$. This model includes the more restrictive models of Jukes-Cantor
(where $b^e=c^e=d^e$, \cite{JC69}) and Kimura 2-parameter model
($b^e=d^e$, \cite{Kimura1980})

In the algebraic geometry setting, the
Kimura 3-parameter model is given by
the polynomial map
\begin{eqnarray*}\prod_{e \in E(T)}\triangle^{3}  & \lra & \triangle^{4^n-1}\\
((a^e,b^e,c^e,d^e))_{e \in E(T)} & \mapsto & (p_{x_1 \dots
x_n})_{x_1, \dots,x_n \in \Sigma^n}
  \end{eqnarray*}
where $p_{x_1, \dots,x_n}$ is given by (\ref{prob}) and $\triangle^d$
denotes the standard $d$-dimensional simplex in $\RR^{d+1}$.
As we are interested in algebraic varieties, instead of
restricting to the simplex, we also consider this polynomial map
as
$$\begin{array}{ccc}\prod_{e \in E(T)}\C^{4}  & \lra & \C^{4^n}\\
((a^e,b^e,c^e,d^e))_{e \in E(T)} & \mapsto & (p_{x_1 \dots
x_n})_{x_1, \dots,x_n \in \Sigma^n}
  \end{array}.$$

One of the goals in
\textit{phylogenetic algebraic
geometry} is determining the ideal of
the closure of the image of this
polynomial map. Knowing the generators
of this ideal provides tools for
phylogenetic inference. See for example
\cite{CFS} and \cite{Eriksson05} where
some of these methods for phylogenetic
inference have been proposed. In order
to find the generators of this ideal,
it is extremely useful to perform
change of coordinates as we explain
below.

\subsection{Fourier transform}\label{subsecfourier}
The models described above are known as
\emph{group-based} models because if
the nucleotides are thought of as the
elements of the group \[H=\ZZ/(2)
\times \ZZ/(2)\] (namely $A=(0,0)$,
$C=(1,0)$, $G=(0,1)$, $T=(1,1)$) then
the entries $\{S^e_{g,h}\}_{g,h\in H}$
in the substitution matrices $S^e$ can
be expressed as functions of the group
$f^e(h-g)$ (see \cite{Sturmfels2005}
for details). For the Kimura
3-parameter model, the function $f^e$
is
$$ f^e(h)=\left\{\begin{array}{lr}a^e & \textrm{ if } h =(0,0) \\
b^e & \textrm{ if } h =(1,0) \\
c^e & \textrm{ if } h =(0,1) \\
d^e & \textrm{ if } h =(1,1) \\
\end{array}\right.$$ As a consequence,
probabilities $p_{x_1  \dots x_n}$ can
also be thought of as functions on $H
\times \dots \times H$. In what
follows, when we add nucleotides, we
mean addition in the group $H$.

One of the main properties of 
group-based models is that a discrete
Fourier transform simplifies the
expression in the probabilities
(\ref{prob}). We briefly recall how
this Fourier transform works  and we
refer to \cite{Sturmfels2005} and
\cite{CGS} for more details. Given a
function $f:H \lra \C$, its discrete
Fourier transform is the function
${\hat{f}:H^{\vee}=Hom{(H,\C^*)} \lra \C}$
defined by
$$\hat{f}(\chi)  \quad = \quad \sum_{g \in H} \chi(g)f(g).$$
The Fourier transform turns convolution into multiplication and
this allows to simplify the expression of joint probabilities:
\begin{thm}[Evans-Speed \cite{Evans1993}]
Let $p(g_1, \dots, g_n)$ be the joint distribution of a
group-based model for a phylogenetic tree $T$, then its Fourier
transform has the form
\begin{equation}\label{fourierparam}
q(\chi_1, \ldots, \chi_n) \quad = \quad
 \!\!\!\prod_{e \in E(T)} \!\!\widehat{f^{e}} \big( \!
\prod_{l  \in \{\textrm{leaves below }
e\}} \chi_l \big)
\end{equation}
\end{thm}
As $H$ and its dual  $H^{\vee}$ are
isomorphic groups we can identify
$(\chi_1, \ldots, \chi_n)$ with the
corresponding tuple $(g_1, \ldots,
g_n)$. From now on, $q(\chi_1, \ldots,
\chi_n)$ will be denoted as $q_{g_1
\dots g_n}$. In the additive notation
of the group $H$, one can  rewrite
expression (\ref{fourierparam}) as
\begin{equation*}
q_{g_1, \ldots, g_n} \quad = \quad
 \!\!\!\prod_{e \in E(T)} \!\!\widehat{f^{e}} \big(m(e) \big)
\end{equation*}
where $m(e)=\sum_{l  \in
\{\textrm{leaves below } e\}}g_l$. The
Fourier transform is a linear
coordinate change given by
\begin{eqnarray}\label{auxfourier}
q_{g_1 \cdots g_n} = \sum_{j_1, \ldots, j_n} \chi^{g_1}(j_1) \cdots \chi^{g_n}(j_n)  p_{j_1 \cdots j_n}
\end{eqnarray}
where $\chi^{i}$ is the character of the group associated to the
$i$th group element:
$$
\begin{array}{r|rrrr}
 & A & C & G & T \\
\hline
\chi^A & 1 & 1  & 1 & 1\\
\chi^C & 1 & -1 & 1 & -1 \\
\chi^G & 1 & 1 & -1 & -1 \\
\chi^T & 1 & -1 & -1 & 1
\end{array}
$$
In this new coordinate system, the Fourier transforms
$\widehat{f^{e}}$ of the substitution functions are the new
parameters of the model. For the Kimura 3-parameter model, these
Fourier transforms become
$$ \widehat{f^e}(h)=\left\{\begin{array}{lr}a^e+b^e+c^e+d^e & \textrm{ if } h =(0,0) \\
a^e-b^e+c^e-d^e & \textrm{ if } h =(1,0) \\
a^e+b^e-c^e-d^e & \textrm{ if } h =(0,1) \\
a^e-b^e-c^e+d^e & \textrm{ if } h =(1,1) \\ \end{array}\right.$$

As before, we think of these
substitution functions as matrices.
Therefore, the parameters in Fourier
coordinates are diagonal matrices
$$P^e=\left(\begin{array}{cccc}
P^e_A & 0 & 0 & 0\\
0 & P^e_C & 0 &0\\
0 & 0 & P^e_G & 0 \\
0 & 0 & 0& P^e_T
\end{array}\right)$$
where $P^e_A= a^e+b^e+c^e+d^e$, $P^e_C=a^e-b^e+c^e-d^e$,
$P^e_G=a^e+b^e-c^e-d^e$, $P^e_T=a^e-b^e-c^e+d^e$.

From now on, $P^e$ will indistinctly denote this diagonal matrix
or the vector $(P^e_A,P^e_C,P^e_G,P^e_T)$ and we will restrict to
Fourier coordinates. It is not difficult to see that if $g_1
+\dots +g_n \neq 0$, then $q_{g_1 \dots g_n}=0$ (cf.
Proposition 29 of \cite{Sturmfels2005}). Therefore, the polynomial map we are
interested in is
\begin{eqnarray}\label{paramfourier}\varphi: \prod_{e \in E(T)}\C^{4}  & \lra &\C^{4^{n-1}}
\\
(P^e_{A}, P^e_{C},P^e_{G},P^e_{T})_e & \mapsto &(q_{x_1 \dots
x_n})_{x_1, \dots, x_n } \nonumber
\end{eqnarray}
where $x_1+ \dots +x_n=0$, $ q_{x_1 \dots x_n}=\prod_{e \in E(T)}P^e_{m(e)}$, and
$m(e_l)=x_l$ if $e_l$ ends at leaf $l$.

\begin{notation}\label{notasimplex} The image of the standard simplex $\triangle^d$
 under the Fourier change of coordinates will be denoted
as $\mathbf{\Delta}^d$.
For a picture of $\mathbf{\Delta}^3$, see figure \ref{simplex}. Notice
that the Fourier change of parameters transforms the hyperplane $a^e+b^e+c^e+d^e=1$ into $P^e_A=1$, for all $e \in E(T)$.
As we will be interested in coordinates $\{q_{x_1,\ldots,x_n}\}_{x_1+\ldots+x_n=0}$, we will focus on the simplex $\D^{4^{n-1}-1}$, which coincides with the projection of the Fourier transform of the simplex $\triangle^{4^{n}-1}$ onto this set of coordinates . As before, note that the Fourier change of coordinates transforms the hyperplane $\sum_{p_{x_1 \dots x_n}}p_{x_1 \dots x_n}=1$ into $q_{A \dots A}=1$. 
\end{notation}

\subsection{Kimura variety}
\label{Kimuravar} In subsections
\ref{subsecbio} and \ref{subsecfourier}
we have assumed that the tree is
rooted. However, in the Kimura
3-parameter model the matrices $S^e$
are symmetric and therefore
parameterization (\ref{paramfourier})
does not depend on the orientation of
the tree or the position of the root.
One can even think of the root as being
one of the leaves and in this case the
root is one of the observed variables.
In what follows we will consider
unrooted trees. This is due to the
issue of \textit{identifiability}  that induces
the use of unrooted trees for
phylogenetic inference. Roughly speaking, the question addressed by identifiability is whether observation data
of character states at the leaves of the tree contain enough information in order to uncover the topology and the parameters of the model (see \cite{Chang96}). This means that there precisely exists one topology and one set of parameters of the model that explain the data. The
identifiability of the Kimura
3-parameter model has been established
by Steel, Hendy and Penny in
\cite{Steel98}.

The reformulation of
the parameterization in Fourier
coordinates for unrooted trees is given
by the lemma \ref{unrooted} below.
Before stating it, we introduce some
notation:
\begin{notation} Given an  interior node
$v\in N(T)$, denote by
$e^1_v,e^2_v,e^3_v$ the three edges
coincident at $v$ (see (a) in figure
\ref{defnxv}). Given three elements of
the group
$x_{e^1_v}$,$x_{e^2_v}$,$x_{e^3_v} \in
H$ associated to these edges, we define
\[x(v):=x_{e^1_v}+x_{e^2_v}+x_{e^3_v}\]
as a sum in $H$.

\begin{figure}
\begin{center}
\psfrag{a}{$e^1_v$}\psfrag{b}{$e^2_v$}
\psfrag{c}{$e^3_v$}\psfrag{v}{$v$}
\psfrag{d}{$e^1_w$}\psfrag{e}{$e^2_w$}
\psfrag{f}{$e^3_w$}\psfrag{w}{$w$}
\psfrag{u}{$T_2$}\psfrag{n}{$T_3$}
\psfrag{(x)}{(a)}\psfrag{(y)}{(b)}
\includegraphics[scale=0.6]{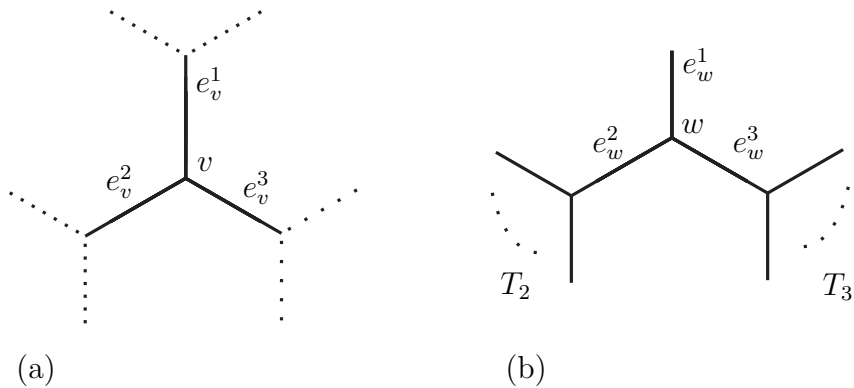}
\end{center}
\caption{\label{defnxv} (a) The three
edges $e^1_v,e^2_v,e^3_v$ coincident at
an interior node $v\in N(T)$. (b) The
two trees obtained from $T$ in the proof of Lemma \ref{unrooted}.}
\end{figure}
\end{notation}

\begin{lem}\label{unrooted}
Let $T$ be an (unrooted) tree with $n$
leaves. Then the parameterization of
the Kimura 3-parameter model in Fourier
coordinates is given by

\begin{equation}\label{qunrooted} q_{x_1 \dots
x_n}=\prod_{e \in E(T), \, x(v)=0 \;
\forall \, v \in N(T)}P^e_{x_e},
\end{equation} where $x_1+ \dots
+x_n=0$ and $x_{e_l}=x_l$ if $e_l$ is
the terminal edge corresponding to the
leaf~$l$.
\end{lem}

\begin{proof}
We already know that the parameterization for rooted trees is
independent of the root placement, so we root the tree $T$ at leaf
$1$. Then we need to prove that
$$\prod_{e \in E(T), x(v)=0 \;
\forall \, v \in N(T)}P^e_{x_e}=\prod_{e \in E(T)}P^e_{m(e)}.$$ In
other words, we want to prove that the condition $x(v)=0$ for all
interior nodes $v \in N(T)$ is equivalent to $x_e=m(e)$ for all
edges $e \in E(T)$.

We first assume that for all interior
nodes $v$, $x(v)=0$. We proceed by
induction on the number of leaves of
the tree. Let $w$ be the only node next
to the leaf $1$ and let $e^1_w$ be the
terminal edge connecting 1 and $w$.
Write $T_2$ and $T_3$ for the two
connected components obtained when
removing $e^1_w$ from $T$ and adding
$W$ as a root (see (b) of figure
\ref{defnxv}). By induction hypotheses
on $T_2$ and $T_3$, we have
$x_{e^2_w}=m(e^2_w)$,
$x_{e^3_w}=m(e^3_w)$, and $x_e=m(e)$
for all edges but $e_1$. It remains to
check that $x_{e_1}=m(e_1)$. This
follows using the hypothesis that
$x(w)=0$. Indeed, $x(w)=0$ implies
$x_1=x_{e^2_w}+x_{e^3_w}$, which in
turn is equal to $m(e^2_w)+m(e^3_w)$
and hence equal to $x_2+ \dots +x_n$.

In order to prove the converse we
assume that $x_e=m(e)$ for all edges $e
\in E(T)$. By induction hypothesis on
the trees $T_2$ and $T_3$, the
condition $x(v)=0$ holds for all
interior nodes in $T$ but $w$. We just
need to show that
$x_{e_1}+x_{e^2_w}+x_{e^3_w}=0$. Our
hypothesis implies that
$x_{e_1}+x_{e^2_w}+x_{e^3_w}=m(e_1)+m(e^2_w)+m(e^3_w)$,
and this vanishes because
$m(e_1)=m(e^2_w)+m(e^3_w)$ by
definition of $m$.
\end{proof}

\begin{rem}
In expression (\ref{qunrooted}), the indices $x_e$ associated to
edge $e$ are completely determined by condition $x(v)=0$, $\forall
v \in N(T)$. Indeed, as at the terminal edges $e_l$ we have
imposed $x_e=x_l$, condition $x(v)=0$ for nodes that join a cherry
to the tree determines the value $x_e$  at those edges that join a
cherry to the tree. Performing the same process from the exterior
of the tree to the interior, one assigns a unique value to every
edge. Condition $x_1+ \dots +x_n=0$ guarantees that this
assignment is consistent at all interior nodes.
\end{rem}

\begin{example} \label{claw}For the unrooted 3-leaf claw tree (see (a) of figure
\ref{clawtree}), the parameterization
$\varphi$ of (\ref{paramfourier}) in
Fourier coordinates is given by
$$q_{x_1x_2x_3}=P^{e_1}_{x_1}P^{e_2}_{x_2}P^{e_3}_{x_3}$$
if $x_1+x_2+x_3=0$.

\begin{figure}
\begin{center}
\psfrag{a}{$e_1$} \psfrag{b}{$e_2$}
\psfrag{c}{$e_3$}
\psfrag{d}{$e_4$}\psfrag{e}{$e$}
\psfrag{(x)}{(a)}\psfrag{(y)}{(b)}
\includegraphics[scale=0.6]{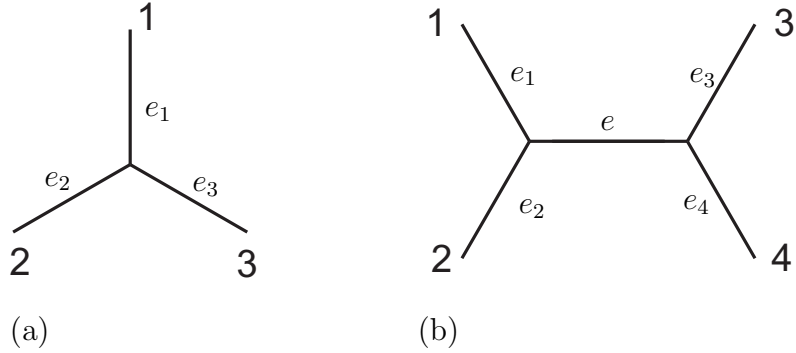}
\end{center}
\caption{\label{clawtree}  (a) The
3-leaf claw tree. (b) The 4-leaves
trivalent unrooted tree.}
\end{figure}
\end{example}

\begin{example}
For the unrooted tree with 4 leaves
(see (b) of figure \ref{clawtree}) the
parameterization in Fourier coordinates
is given by
$$q_{x_1x_2x_3x_4}=P^{e_1}_{x_1}P^{e_2}_{x_2}P^e_{x_1+x_2}P^{e_3}_{x_3}P^{e_4}_{x_4}$$
if $x_1+x_2+x_3+x_4=0$.
\end{example}

\begin{notation}
From now on, we write $V$ for the
closure of the image of
\begin{eqnarray*}\varphi: \prod_{e \in E(T)}\C^{4}  & \lra &\C^{4^{n-1}}
\\
(P^e_{A}, P^e_{C},P^e_{G},P^e_{T})_e &
\mapsto &(q_{x_1 \dots x_n})_{x_1+
\dots +x_n =0}
\end{eqnarray*} where
\begin{equation*}
q_{x_1 \dots x_n}=\prod_{e \in E(T), \,
x(v)=0 \; \forall \, v \in
N(T)}P^e_{x_e}, \end{equation*}  and
$x_{e_l}=x_l$ if $e_l$ is the terminal
edge corresponding to leaf $l$. We
will denote $\Q_n$ the following set of
indeterminates
\[\Q_n= \{q_{x_1\dots x_n} |
{x_1+\dots+x_n= 0}\textrm{ in }H\}.\]
\end{notation}

\vspace{4mm}

The affine coordinate ring $A(V)$ of $V$ is isomorphic to the
$\C$-algebra
$$\C[\{\prod_{
\begin{tiny}
\begin{array} {cc} e \in E(T) \\
x(v)=0 , \, \forall  v \in N(T)
\end{array}
\end{tiny}
}P^e_{x_e} \mid x_{e_l}=x_l \; \forall \, l\in N(T)  \}_{\, x_1 \dots x_n \in \Sigma}],$$
because $A(V)$ is defined as the image
of  the  morphism of $\C$-algebras
\begin{eqnarray}\label{morph}
\theta:\C[\Q_n] & \lra & \C[\{P^e_{x} \mid {e \in E(T) , x \in \Sigma}\}]\\
q_{x_1\dots x_n}& \mapsto & \prod_{e
\in E(T), \, x(v)=0 \; \forall \, v \in
N(T)}P^e_{x_e} \nonumber
\end{eqnarray}
where $x_{e_l}=x_l$ for all $l\in
N(T).$ The toric ideal defining $V$ is
the kernel of this morphism and we
denote it as $I_V$. Sturmfels and
Sullivant gave an algorithm to
construct a  set of generators of $I_V$
in \cite{Sturmfels2005}. We note that
as the map $\theta$ is homogeneous, $V$
is actually a cone over a projective
variety.

The variety we are interested in is $V
\cap \{q_{A\dots A}=1\}$ because the
simplex in Fourier coordinates is
contained in the hyperplane $q_{A\dots
A}=1$.
\begin{defn}\rm The \textit{Kimura
variety} of the phylogenetic tree $T$ is $$W:=V \cap \{q_{A\dots
A}=1\}\,.$$
\end{defn}
\begin{lem}\label{Kimvar1}
The Kimura variety $W$ is the closure
of the parameterization $\varphi$
restricted to $\prod_{e \in E(T)
}(\C^{4}\cap \{P^e_A=1\} )$.
\end{lem}

\begin{proof}
Write $\varphi_1$ for this restriction,
and let $\rho$ be the morphism of
$\C$-algebras defined as:
$$\begin{array}{ccc} \rho:\C[\Q_n] & \lra & \C[\Q_n \setminus\{q_{A\dots A}\}] \\
f(q_{A \dots A}, q_{A \dots CC}, \dots) & \mapsto &f(1, q_{A \dots
CC}, \dots)\end{array}
$$
Then the affine coordinate ring of $W$
is \[A(W) = \C[\Q_n] / (I_V+(q_{A \dots
A}-1)) \simeq \C[\Q_n \setminus \{q_{A
\dots A}\}]/(\rho(I_V).\]On the other
hand, note that $\varphi_1$ induces a
morphism of $\C$-algebras
\begin{eqnarray*}
\theta_1:\C[\Q_n \setminus \{q_{A
\dots A}\}]  \lra \C[\{P^e_{x} \mid {e \in E(T)} \textrm{
and } x \in \Sigma\setminus \{A\}\}]
\end{eqnarray*}
so that the following diagram commutes
\begin{equation}\label{diagram}\begin{array}{ccccccc}
0 & \lra & I_V & \lra & \C[\Q_n] & \srel{\theta}{\lra} & \C[\{P^e_{x} \mid {e \in E(T)} \mbox{ and } x \in \Sigma\}] \\
 &      & \downarrow \rho& &  \downarrow \rho & & \downarrow \rho' \\
0 & \lra & \rho(I_V) & \lra & \C[\Q_n \setminus\{q_{A\dots A}\}] &
\srel{\theta_1}{\lra} & \C[\{P^e_{x} \mid {e \in E(T)} \textrm{
and } x \in \Sigma\setminus \{A\}\}]
 \end{array} \,.
\end{equation}
Here $\rho'$ sends $P^e_A$ to $1$ for each $e \in E(T)$ and is the identity on the other indeterminates.
Write $X\subset \C^{4^n}$ for the
closure of the image of $\varphi_1$.
Then, the affine coordinate ring of $X$
is $A(X)=\C[\Q_n]/Ker (\rho'\circ
\theta)$. Since the above diagram
commutes, we have that \[\rho'\circ
\theta=\theta_1\circ
\rho=\rho^{-1}(\rho(I_V))=I_V+(q_{A\ldots
A}-1)\]and $X=W$.
\end{proof}

\subsection{Biologically meaningful
points}\label{biology}
\begin{figure}
\begin{center}
\includegraphics[scale=0.6]{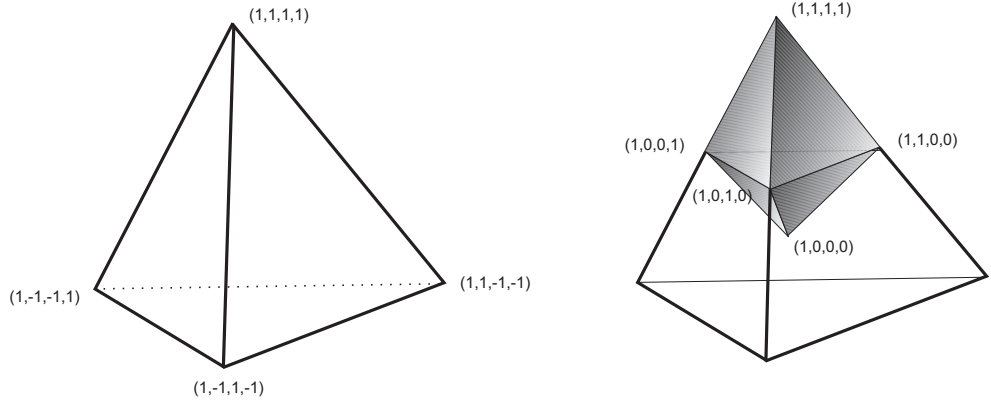}
\end{center}
\caption{\label{simplex}  The polytope
$\D^{3}_0$ in the simplex $\D^{3}$. It
represents the points of $\D^{3}$
having non-negative Fourier coordinates.}
\end{figure}

Here, we introduce some notation and we
give a biological interpretation on the
points of some polytopes in the
simplices defined in subsection
\ref{subsecfourier}. As above, let $T$
be a phylogenetic tree with $n$ leaves
and write $E(T)$ for the set of its
edges. For any $e\in E(T)$, let
$\D^{3}$ be the parameter simplex
associated to it (see Notation (\ref{notasimplex})). Write
\[\D^{3}_0=\{P=(1,P_C,P_G,P_T)\in \D^{3} \mid P_C,P_G,P_T\geq 0\}\]for
the set of non-negative points in
$\D^3$ i.e. the polytope delimited
by the vertex $(1,1,1,1)$, the points
$(1,1,0,0),(1,0,1,0),(1,0,0,1)$ and the
centroid $(1,0,0,0)$ of $\D^3$ (see figure \ref{simplex}), all
of them  written in Fourier coordinates. We write
also $\D^{3}_+=\{P\in \D^{3} \mid
P^e_C,P^e_G,P^e_T> 0\}$. Any point of
$\D^{3}_0$ represents a substitution
matrix $S^e$ satisfying that the
probability of no mutation ($a^e$) plus
the probability of any
mutation in the site ($b^e$, $c^e$ or
$d^e$) is bigger or equal than $1/2$. In
particular, the probability of no
mutation is bigger or equal than the
probability of any particular mutation.
Since this is a reasonable hypothesis
if we work with realistic data, we will
call the points in $\prod_{e\in
E(T)}\D^{3}_+$ the
\textit{biological meaningful parameters
of the model}. We also write
\begin{eqnarray*}
 \varphi_+:\prod_{e\in E(T)}\D^{3}_+ & \lra & W \\
(P^e_{A}, P^e_{C},P^e_{G},P^e_{T})_e & \mapsto &(q_{x_1 \dots
x_n})_{x_1, \dots, x_n }
\end{eqnarray*}
for the restriction of $\varphi$ in
(\ref{paramfourier}) to these
parameters.  Its image is contained  in
$W_+:=W\cap \D_+$, where $\D_+=\{q\in
\D^{4^{n-1}-1} \mid q_{x_1 \ldots x_n}
> 0,\, \sum_i x_i=0\}$ is the set of points with positive
coordinates of the polytope delimited
by the point
$\mathbf{1}_n=(1,\ldots,1)$ (which is
one of the vertices on the simplex
${\D^{4^{n-1}-1}}$), the points
$q_i=(1,e_i)_{i=1,\ldots,4^{n-1}-1}$
where $e_i=(0,\ldots,1,\ldots,0)\in
\C^{4^{n-1}-1}$ and the centroid
$(1,0,\ldots,0)$ of $\D^{4^{n-1}-1}$.
We call
the points of $W_+$ the
\textit{biologically meaningful points
of the model}. This name will be justified in forthcoming Remark \ref{preim}, where we show that $W_+$ equals the image of $\varphi_+$.

\section{The geometry of Kimura 3-parameter model}
Let $V \subset \C^{4^n}$ be the affine
variety associated to a tree $T$ of $n$
leaves, as defined above. In this
section we are going to determine the
singular points of the Kimura variety
$W=V \cap \{q_{A\dots A}=1\}$. To this
aim, we will first prove that $V$ is
isomorphic to the quotient of
$(\C^{4})^{2n-3}$ by the action of a
certain abelian group.

In order to simplify notation, we
briefly recall  the notion of
\textit{multigrading} and refer to
Chapter 8 of \cite{MillerSturmfels} for
a nice introduction to multigraded
polynomial rings (we also refer to
\cite{MillerSturmfels} for the
correspondence between toric varieties
and quotients, although we need little
preliminary knowledge in this subject
for the results of this section.)

\begin{notation}\label{multigrading}
Let $M$ be a monomial in
$S=\C[\{P^e_{x} \mid {e \in E(T)}
\textrm{ and } x \in \Sigma\}]$. Then,
$M$ has the form $M= \prod_{e \in E(T)}
(P^e)^{\mathbf{i}(e)}$ where each
$\mathbf{i}(e)=(i(e)_A,i(e)_C,i(e)_G,i(e)_T)$
is composed of natural numbers. This
notation means that  \[M=\prod_{e \in
E(T), x \in \Sigma }
(P^e_x)^{i(e)_x}.\] We call
$\deg(\mathbf{i}(e))=i(e)_A+i(e)_C+i(e)_G+i(e)_T$.
Each indeterminate in $S$ has a natural
multidegree in $\ZZ^{|E(T)|}$ defined
as
$$\deg(P^e_x)=(0, \dots,0,\srel{e)}{1},0, \dots, 0) \quad \mbox{for any }x\in \Sigma.$$
Given a monomial $M \in S$ as above, we
call $\deg(M):=(\deg(\mathbf{i}(e)))_{e
\in E(T)}$. Note that the image of
$\theta$ of (\ref{morph}) is generated
by monomials of degree $d \cdot(1,
\dots, 1)$, so that they are
multi-homogeneous with respect to the
given grading.
\end{notation}

\begin{notation} From now on, $\ZZ/(2)$ means additive group whereas $\ZZ_2$
means multiplicative group.
\end{notation}

\subsection{The 3-leaves case}
We start by studying the  case $n=3$
(see example \ref{claw}). We call $V_3$
the corresponding affine variety in
$\C^{16}$. The parameterization
$\varphi$ in this case is:
$$ \begin{array}{ccc} \varphi:\C^{12} & \lra & \C^{16}\\
\scriptstyle{((P^{e_1}_A,P^{e_1}_C,P^{e_1}_G,P^{e_1}_T),(P^{e_2}_A,P^{e_2}_C,P^{e_2}_G,P^{e_2}_T),(P^{e_3}_A,P^{e_3}_C,P^{e_3}_G,P^{e_3}_T))} &\mapsto
& \scriptstyle{(P^{e_1}_xP^{e_2}_yP^{e_3}_z)_{\{x+y+z =0 \mid x,y,z \in
\Sigma\}}}\end{array}$$

In the next result we prove that $V_3$ is an affine GIT quotient
\cite[chapter 10]{MillerSturmfels}.
\begin{prop}\label{prop3fulles}
$V_3$ is isomorphic to the affine GIT quotient $\C^{12} // G$
where the group
\[G=\{(\lambda_1, \lambda_2, \lambda_3, \eps, \delta) \mid \lambda_i \in
\C^*, (\eps, \delta) \in \ZZ_2 \times \ZZ_2, \lambda_1\lambda_2
\lambda_3=1 \}\simeq (\C^*)^2 \times \ZZ_2 \times \ZZ_2\] acts on
$\C^{12}$ sending $(P^{e_1},P^{e_2},P^{e_3})$ to \[(\lambda_1(P^{e_1}_A,\eps
P^{e_1}_C,\delta P^{e_1}_G,\eps \delta P^{e_1}_T),\lambda_2(P^{e_2}_A,\eps
P^{e_2}_C,\delta P^{e_2}_G,\eps \delta P^{e_2}_T),\lambda_3(P^{e_3}_A,\eps
P^{e_3}_C,\delta P^{e_3}_G,\eps \delta P^{e_3}_T)).\]
\end{prop}

In order to prove this proposition, we need a technical lemma that
we state separately for future reference.

\begin{lem}\label{lemaindex}
Let $\mathbf{i}=(i_A,i_C,i_G,i_T)$,
$\mathbf{j}=(j_A,j_C,j_G,j_T)$, $\mathbf{k}=(k_A,k_C,k_G,k_T)$ be
4-tuples in $\NN^4$. Then the set of indices
$(\mathbf{i},\mathbf{j},\mathbf{k})$  that satisfy
\begin{equation}\label{eqs}
\left\{\begin{array}{c} i_A+i_C+i_G +i_T=1\\
j_A+j_C+j_G +j_T=1 \\ k_A+k_C+k_G +k_T=1 \\
i_C+i_T+j_C+j_T+k_C+k_T=0 \textrm{ in } \ZZ/(2)  \\
i_G+i_T+j_G+j_T+k_G+k_T =0 \textrm{ in } \ZZ/(2)
\end{array}\right.\end{equation} is equal to
$$\{(\mathbf{i},\mathbf{j},\mathbf{k}) \, \mid \,i_x=1, j_y=1,k_z=1,
 \deg(i)=\deg(j)=\deg(k)=1, x+y+z=0 \textrm{ in } H\}.$$

\end{lem}

\begin{proof}
Let $(\mathbf{i},\mathbf{j},\mathbf{k})$ satisfy (\ref{eqs}). The
first three equations imply that for each index, there is just one
letter in $\Sigma$ such that the corresponding entry is non-zero,
and in fact, equal to one. We write $i_x=1$, $j_y=1$, $k_z=1$. As
in the first section we think of the letters in $\Sigma$ as
elements in $\ZZ/(2) \times \ZZ/(2)$: $A=(0,0)$, $C=(1,0)$,
$G=(0,1)$, $T=(1,1)$. We call
\[\begin{array}{lll}
 I_{AC}=i_A+i_C, &  I_{CT}=i_C+i_T, &  I_{GT}=i_G+i_T, \\
 J_{AC}=j_A+j_C, &  J_{CT}=j_C+j_T, &  J_{GT}=j_G+j_T, \\
 K_{AC}=k_A+k_C, &  K_{CT}=k_C+k_T, &  K_{GT}=k_G+k_T.
\end{array}\]
In this setting, $i_x=1$ if and only if $x=(I_{CT},I_{GT})$ in
$\ZZ/(2) \times \ZZ/(2).$ Similarly, $j_y=1$ (resp. $k_z=1$) if
and only if $y=(J_{CT},J_{GT})$ (resp. $z=(K_{CT},K_{GT})$  ) in
$\ZZ/(2) \times \ZZ/(2)$. The last two equations in (\ref{eqs})
can be written as
\begin{equation*}
\left\{\begin{array}{c} I_{CT}+J_{CT}+K_{CT}=0 \textrm{ in } \ZZ/(2)  \\
I_{GT}+J_{GT}+K_{GT}=0 \textrm{ in } \ZZ/(2)
\end{array}\right. \end{equation*}
They imply that $x+y+z=0$ in $\ZZ/(2) \times \ZZ/(2)$.

As for the other inclusion, if
$(\mathbf{i},\mathbf{j},\mathbf{k})$ are three 4-tuples in $\NN^4$
of degree 1 ($\deg(\mathbf{i}) = \deg(\mathbf{j}) =
\deg(\mathbf{k}) = 1$) whose non-vanishing indices
$i_x=1$,$j_y=1$, $j_z=1$ satisfy $x+y+z=0$. Then the first three
equations in (\ref{eqs}) are clearly satisfied. The last two
equations also hold because $(i_C+i_T,i_G+i_{T})=x$,
$(j_C+j_T,j_G+j_{T})=y$, $(k_C+k_T,k_G+k_{T})=z$ and $x+y+z=0$ in
$\ZZ/(2) \times \ZZ/(2)$.
\end{proof}

We prove the proposition above.

\vspace{3mm}
\noindent \textsc{Proof of Proposition \ref{prop3fulles}.}
Recall that $\C^{12} //G$ is defined as
$Spec(S^G)$, the spectrum of the ring of invariants $S^G$ where
$S=\C[P^{e_1}_A, \dots, P^{e_1}_T, P^{e_2}_A, \dots, P^{e_2}_T,  P^{e_3}_A, \dots
P^{e_3}_T]$. This ring is generated as a $\C$-algebra by those
monomials invariant by the action of $G$.

Monomials in $S$ are of the form
$(P^{e_1})^\mathbf{i}(P^{e_2})^\mathbf{j}(P^{e_3})^\mathbf{k}$ where
$\mathbf{i},\mathbf{j},\mathbf{k}$ are sets of natural numbers
$\mathbf{i}=(i_A, i_C,i_G,i_T)$, $\mathbf{j}=(j_A, j_C,j_G,j_T)$,
$\mathbf{k}=(k_A, k_C,k_G,k_T)$ and $(P^{e_1})^\mathbf{i}$ means
$(P^{e_1}_A)^{i_A}(P^{e_1}_C)^{i_C}(P^{e_1}_G)^{i_G}(P^{e_1}_T)^{i_T}$. A monomial
is invariant under the action of $G$ if and only if for any
$(\lambda_1, \lambda_2, \lambda_3, \sigma=(\eps,\delta))$ in $G$
we have
$$\lambda_1^{i_A+\dots +i_T}\lambda_2^{j_A+\dots+ j_T}\lambda_3^{k_A+\dots +k_T}\eps^{i_C+i_T+j_C+j_T+k_C+k_T}\delta^{i_G+i_T+j_G+j_T+k_G+k_T}=1$$
This happens if and only if
$$\left\{
\begin{array}{c} i_A+\dots+ i_T=j_A+\dots +j_T=k_A+\dots +k_T \\
i_C+i_T+j_C+j_T+k_C+k_T=0 \textrm{ in } \ZZ_2  \\
i_G+i_T+j_G+j_T+k_G+k_T =0 \textrm{ in } \ZZ_2\end{array}\right.$$
Therefore $S^G$ is minimally generated as a $\C$-algebra by those
monomials $(P^{e_1})^\mathbf{i}(P^{e_2})^\mathbf{j}(P^{e_3})^\mathbf{k}$ that
satisfy
$$\left\{\begin{array}{c} i_A+i_C+i_G +i_T=1\\
j_A+j_C+j_G +j_T=1 \\ k_A+k_C+k_G +k_T=1 \\
i_C+i_T+j_C+j_T+k_C+k_T=0 \textrm{ in } \ZZ_2  \\
i_G+i_T+j_G+j_T+k_G+k_T =0 \textrm{ in } \ZZ_2\end{array}\right.$$
By the lemma following this proof, this set of monomials is
precisely
$$\{P^{e_1}_xP^{e_2}_yP^{e_3}_z \mid x+y+z=0 \textrm{ in } \ZZ/(2) \times \ZZ/(2)\}.$$
Therefore $S^G$ is the finitely generated $\C$-algebra
$$\C[\{P^{e_1}_xP^{e_2}_yP^{e_3}_z \mid x+y+z=0 \textrm{ in } \ZZ/(2) \times
\ZZ/(2)\}]$$ which is isomorphic to the affine coordinate ring of
$V$.
\qed

\subsection{The general case}
Now we generalize Proposition
\ref{prop3fulles} to trees with an
arbitrary number $n$ of leaves. Recall
that the number of edges in such a tree
is $2n-1$. By a \emph{path} $\sigma$ in
$T$ we mean a minimal subgraph of $T$
connecting two nodes (interior nodes or
leaves). We write $\sigma=\{s_1,\ldots,
s_r\}$ for the sequence of edges in
$\sigma$. Let
\[G\subset (\C^*\times \ZZ_2 \times
\ZZ_2)^{2n-3}\]be defined as the
subset composed of the elements
$(\lambda_e,\eps_e,\delta_e)_{e \in
E(T)}$ such that $\prod_{e \in E(T)}
\lambda_e=1$ and satisfying the
following condition
\begin{itemize}  \item[($\ast$)] for  any path $\sigma=\{s_1,\dots, s_r\}$ between two leaves in $T$,
$\prod_{i=1}^r \eps_{s_i}=1$ and $
\prod_{i=1}^r \delta_{s_i}=1$.
\end{itemize}
It is immediate to see that  the
natural product induces a group
structure in $G$. Moreover, we have

\begin{lem}
The group $G$ is  isomorphic to
$(\C^*)^{2n-4}\times (\ZZ_2 \times
\ZZ_2)^{n-2}$.
\end{lem}

\begin{proof}
For any interior node $v$ of $T$, let $e_v^1,e_v^2,e_v^3$ be the
edges incident at $v$, and write
$\eps(v)=\{\eps_e,\delta_e\}_{e\in E(T)}$ by taking
$\eps_{e_v^1}=\eps_{e_v^1}=\eps_{e_v^1}=-1$, $\eps_e=1$ for the
remaining edges. Similarly, we define $\delta(v)$.
Let $e_0\in E(T)$ and take the ring homomorphism
\[\psi: (\C^*)^{|E(T)|-1}\times
(\ZZ_2 \times \ZZ_2)^{|N(T)|} \lra (\C^*\times \ZZ_2 \times
\ZZ_2)^{2n-3}\]defined by mapping $((\mu_e)_{e\neq e_0},
(\eps(v),\delta_v)_{v\in N(T)})$ to
$(\lambda_e,\eps_e,\delta_e)_{e\in E(T)}$, where $\lambda_e=\mu_e$
if $e\neq e_0$, $\lambda_{e_0}=(\prod_{e\neq e_0} \lambda_e)^{-1}$
and $\eps_e=\prod_{v\in e}\eps(v)$, $\delta_e=\prod_{v\in
e}\delta(v)$. The image of $\psi$ is $G$ and it is easy to check
that $\psi$ is a monomorphism. The claim follows.
\end{proof}

The main result of this section is the following theorem.

\begin{thm}\label{quotientG} Let $T$  be a tree with $n$ leaves and let $G$ be the group defined
above. Let $G$ act on $\prod_{e\in E(T)}\C^{4}$ by sending
$(P^{e}_A,P^{e}_C,P^{e}_G,P^{e}_T)_{e\in E(T)}$ to
\[(\lambda_e(P^e_A,\eps_e P^e_C,\delta_e P^e_G,\eps_e
\delta_e P^e_T))_{e \in E(T)}.\] Then $V$ is isomorphic to
$(\C^{4})^{2n-3}// G$.
\end{thm}
\begin{proof} We need to  check that the affine coordinate rings of
$V$ and $(\C^{4})^{2n-4}//G$ are isomorphic. If $S$ is the algebra
$S=\C[\{P^e_{x} \mid {e \in E(T)} \textrm{ and } x \in \Sigma\}]$,
we need to check that the ring of invariants $S^G$ is isomorphic
to
 \[\C[\{\prod_{e \in E(T), x(v)=0 \; \forall
 \, v \in N(T)}P^e_{x_e} \, \mid \,  x_l
 \in \Sigma, x_{e_l}=x_l \; \forall \,
 l\in \{1, \dots,n\} \,].\]
Let $M\in S$ be a monomial. Then, $M$ has the form \[M= \prod_{e
\in E(T)} (P^e)^{\mathbf{i}(e)},\]with the notation introduced in
\ref{multigrading}.
This monomial is invariant by the action of $G$ if and only if for
any $(\lambda_e,\eps_e, \delta_e)_{e \in E(T)} \in G$ we have
\begin{equation}\label{eqlambda}1=\prod_{e\in
E(T)}\lambda_e^{\deg(\mathbf{i}(e))}\prod_{e\in
E(T)}\eps_e^{i(e)_C+i(e)_T}\prod_{e\in
E(T)}\delta_e^{i(e)_G+i(e)_T}.\end{equation} As  $\prod_{e\in
E(T)}\lambda_e=1$, equation (\ref{eqlambda}) implies
$\deg(\mathbf{i}(e))=\deg(\mathbf{i}(e'))$ for all $e,e' \in E(T)$
(in the language  of the previous section, this means that $M$ is
multi-homogeneous).
Therefore the  algebra $S^G$ is generated by monomials that
satisfy $\deg(\mathbf{i}(e))=1$ for all edges $e \in E(T)$. We
assume from now on that $M$ satisfies this condition.

Let $v$ be an interior node of $T$, and let $e_v^1,e_v^2,e_v^3$ be
the edges incident at $v$. If we take
$\eps_{e_v^1}=\eps_{e_v^1}=\eps_{e_v^1}=-1$ (resp.
$\delta_{e_v^1}=\delta_{e_v^1}=\delta_{e_v^1}=-1$) and $\eps_e=1$
(resp. $\delta_e=1$) for the remaining edges, condition ($\ast$)
is satisfied. For this particular choice, equation
(\ref{eqlambda}) implies
 $$\left\{
\begin{array}{c}
i(e_v^1)_C+i(e_v^1)_T+i(e_v^2)_C+i(e_v^2)_T+i(e_v^3)_C+i(e_v^3)_T=0 \textrm{ in } \ZZ/(2)  \\
i(e_v^1)_G+i(e_v^1)_T+i(e_v^2)_G+i(e_v^2)_T+i(e_v^3)_G+i(e_v^3)_T=0
\textrm{ in } \ZZ/(2)
\end{array}\right.$$
Lemma \ref{lemaindex} tells us that $M=\prod_{e \in E(T), x(v)=0
\; \forall \, v \in N(T)}P^e_{x_e} $.

We need to prove the converse: if $M=\prod_{e \in E(T), x(v)=0 \;
\forall \, v \in N(T)}P^e_{x_e}$ for some given $x_1, \dots,x_n$,
we shall check that it is invariant  by the action of $G$. In
other words, if $\{\mathbf{i_e}\}_{e \in E(T)}$ denotes the set of
exponents in $M$, we are going to check that equation
($\ref{eqlambda}$) holds. By lemma \ref{lemaindex}, condition
$x_{e^1_v}+x_{e^2_v}+x_{e^3_v}=0$ is equivalent to
\begin{equation}\label{condi}
\begin{array}{c}
i(e_v^1)_C+i(e_v^1)_T+i(e_v^2)_C+i(e_v^2)_T+i(e_v^3)_C+i(e_v^3)_T=0 \textrm{ in } \ZZ/(2)  \\
i(e_v^1)_G+i(e_v^1)_T+i(e_v^2)_G+i(e_v^2)_T+i(e_v^3)_G+i(e_v^3)_T=0
\textrm{ in } \ZZ/(2)
\end{array}
\end{equation}

\vspace{3mm} \noindent \textsc{Claim}: If condition (\ref{condi})
holds for any $v \in N(T)$, the sets $\gamma_{CT}=\{e \in E(T)
\mid i(e)_C+i(e)_T=1\}$ and $\gamma_{GT}=\{e \in E(T) \mid
i(e)_G+i(e)_T=1\}$ are unions of disjoint paths between leaves of
$T$.

\noindent \emph{Proof}: If $x_i=A$ for all leaves, then the set of
conditions $x_{e^1_v}+x_{e^2_v}+x_{e^3_v}=0$ lead to $x_e=A$ for
all $e \in E(T)$, so there is nothing to prove in this case. We
assume that there is a leaf $l$ such that $x_l \neq A$ and we
assume that $x_l\in \{C,T\}$ (if $x_l\in \{G,T\}$ we proceed
analogously). Then $e_l$ belongs to $\gamma_{CT}$.

Let $v$ be the interior node connecting the edge $e_l$ to the rest
of the tree. As $e_l$ is one of the edges intersecting at $v$,
condition
$i(e_v^1)_C+i(e_v^1)_T+i(e_v^2)_C+i(e_v^2)_T+i(e_v^3)_C+i(e_v^3)_T=0
\textrm{ in } \ZZ/(2)$ implies that one of the other two edges
emerging from $v$ also belongs to $\gamma_{CT}$. We call this edge
$e_vw$ and $w$ is the other extreme of the edge. Then for $w$
condition
$i(e_w^1)_C+i(e_w^1)_T+i(e_w^2)_C+i(e_w^2)_T+i(e_w^3)_C+i(e_w^3)_T=0
\textrm{ in } \ZZ/(2)$ again implies that one of the other two
edges emerging from $w$ belongs to $\gamma_{CT}$. We repeat this
process until we end at another leaf of $T$. We note that any two
paths obtained this way are disjoint because an interior node
cannot have three edges in $\gamma_{CT}$. Therefore the claim is
proved.

The claim immediately implies that the
monomial $M$ satisfies equation
($\ref{eqlambda}$), since elements in
the group $G$ are defined by the
condition ($\ast$).
\end{proof}

\begin{cor}\label{finit_group}
The coordinate ring of the Kimura variety $W=V \cap \{q_{A \dots
A}=1\}$ is isomorphic to $S'^{G'}$ where
$$S'=\C[\{P^e_{x} \mid {e \in E(T)} \textrm{ and } x \in \Sigma\setminus \{A\}\}]$$
and $G'= (\ZZ_2 \times \ZZ_2)^{n-2}$ acts as a subgroup of the
group $G$ defined in Theorem \ref{quotientG}. Equivalently, $W$ is
the affine GIT quotient of $\prod_{e \in E(T) }(\C^{4}\cap
\{P^e_A=1\} )$ modulo $G'$,
$$\prod_{e \in E(T) }(\C^{4}\cap
\{P^e_A=1\} )// G'\,.$$
\end{cor}
\begin{proof}
Using diagram (\ref{diagram}) we see that $A(W)$ is isomorphic to
$\rho(A(V))$. By Theorem \ref{quotientG}, we know that $A(V)\simeq
S^G$ and it is enough to prove $\rho(S^G)=S'^{G'}.$

We first prove $\rho(S^G) \subseteq S'^{G'}$. Let
$M((P^e_A,P^e_C,P^e_G,P^e_T)_{e \in E(T)})$ be a monomial in
$S^G$. Then $\rho(M)=M((1,P^e_C,P^e_G,P^e_T)_{e \in E(T)})$. Let
$g'=(\eps_e,\delta_e)_{e \in E(T)}$ be an element in $G'$. We have
that the action of $g$ in $\rho(M)$ is $g \cdot
\rho(M)=M((1,\eps_e P^e_C,\delta_e P^e_G,\eps_e \delta_e P^e_T)_{e
\in E(T)})$.
Take $g=(1,\eps_e,\delta_e)_{e\in E(T)}$, which is an element of
$G$, and notice that
\[g\cdot M=g'\cdot \rho(M).\]
As $M$ is invariant by the action of $G$, it is also invariant by
this element $g'$. Hence $g\cdot \rho(M)=\rho(M).$

In order to prove the other  inclusion we will use the
multigrading notation introduced in \ref{multigrading}. Let
$M((P^e_C,P^e_G,P^e_T)_{e \in E(T)})=\prod_{e \in E(T)}
(P^e)^{\mathbf{i}(e)}$ be a  monomial in $S'^{G'}$. As $S'$ is a
subring of $S$, there is also a multidegree associated to $M$,
namely $\deg(M)=(\deg(\mathbf{i}(e)))_{e \in E(T)}$. Now we make
$M$ multi-homogeneous: let $D$ be equal to $\max_{e \in
E(T)}\deg(\mathbf{i}(e))$ and consider the monomial $N:=\prod_{e
\in E(T)} (P_A^e)^{d -\deg(\mathbf{i}(e))}M$. Then $N$ is a
monomial in $S$ invariant by the action of $G$ because $M$ was
invariant by $G'$ and $N$ is multi-homogeneous (so that equation
(\ref{eqlambda}) holds). Moreover $\rho(N)=M$ and we are done.
\end{proof}

\begin{rem}\label{reflection}
It is worth pointing out that the action of $\mathbb{Z}_2 \times \mathbb{Z}_2$ on $\D^{3}$:
\[(g,P^e)=((\eps,\delta),(1,P_C^e,P_G^e,P_T^e)\mapsto
(1,\eps P_C^e,\delta P_G^e, \eps \delta
P_T^e).\]is just the
reflection relative to some of the axis
going through the centroid of
$\D^{3}$. Namely, the actions of
$g_1=(-1,1)$, $g_2=(1,-1)$ and
$g_3=(-1,-1)$ are the reflections
relative to the $P_G$-axis, the
$P_C$-axis and the $P_T$-axis,
respectively. Thus, if we write
$\D_{x,y}=\{P \in \D^{3} \mid
P^e_x,P^e_y\leq 0\}$ for any $x,y\in
\Sigma$, then
$g_1(\D^{3}_0)=\D_{C,T}$,
$g_2(\D^{3}_0)=\D_{G,T}$ and
$g_3(\D^{3}_0)=\D_{C,G}$.
\end{rem}

\begin{cor}\label{surjective}
The Kimura variety $W$ is the geometric
quotient $$\prod_{e \in E(T)
}(\C^{4}\cap \{P^e_A=1\} )/ G'$$ and
coincides  with the image of
$\varphi_1$ (cf. Lemma \ref{Kimvar1})
\end{cor}

\begin{proof}
By Corollary \ref{finit_group}, the
variety $W$ is the categorical quotient
defined by $Spec(S')^{G'}$. As $G'$ is
a finite group, the orbits of $G'$ are
closed and this categorial quotient is
precisely the geometric quotient
$\prod_{e \in E(T) }(\C^{4}\cap
\{P^e_A=1\} )/G'$ and therefore, it
coincides with the image of $\varphi_1$
(see Example 6.1 of \cite{Dolgachev}).
\end{proof}

From this, we deduce the identifiability of the model (see subsection \ref{biology} and \cite{Chang96}). In particular, we have:

\begin{cor}\label{dim_codim} The Kimura variety $W \subset \CC^{4^{n-1}}$

has dimension $3(2n-3)$ and codimension $4^{n-1}-6n+9$.
\end{cor}

\begin{cor}
Let $q$ be a point in the Kimura variety $W$. Then
\[|\varphi^{-1}(q)|\leq 4^{n-2}\]and the equality holds for generic points.
The same holds for $q\in
W^{\mathbb{R}}=\varphi((\mathbb{R}^3)^{2n-3})$ or $q\in
W_{\mathbf{\Delta}}=\varphi((\mathbf{\Delta}^3)^{2n-3})$.
\end{cor}

\begin{proof}
By Corollary \ref{surjective}, $W$ is the image of $\varphi$ in the commutative diagram
\[\xymatrix{  {\prod_{e \in
E(T) }(\C^{4}\cap \{P^e_A=1\} )} \ar[r]^{\qquad \qquad \varphi}\ar[rd]_{\pi} & {W} \ar[d]^{\cong} \\ &  {\prod_{e \in
E(T) }(\C^{4}\cap \{P^e_A=1\} )/G'}}\]
It follows that if $q\in W$, $\varphi^{-1}(q)$ consists of one
point $(1,P^e_C,P^e_G,P^e_T)_{e\in E(T)}$ and its images under the
action of $G'$. Since $|G'|=4^{n-2}$, the claim follows. Moreover, the image of a point in $\mathbb{R}^3$ (resp. $\D^3$) under the action of $G'$ stays in $\mathbb{R}^3$ (resp. $\D^3$). 
\end{proof}

\begin{rem} \label{preim}
Notice that the pre-images by $\varphi_1$
of the point $\mathbf{1}_n=(1,\ldots,
1)$ are \[\varphi^{-1}(p)=\{(1,\eps_e
P^e_C,\delta_e P^e_G,\eps_e \delta_e
P^e_T)_{e\in E(T)}\mid
(\eps_e,\delta_e)_e \in G'\}.\]Among
all of them, there is only one with
biological interest: $((1,1,1,1)_{e\in
E(T)})$, which represents the situation
where, in probability parameters, the
transition matrices of all edges are
equal to the identity (no mutation
occurs). 
In general, for any point $q\in W_+$ with real coordinates, there is one just preimage of biological interest.

Keeping the notation of subsection
\ref{biology}, it follows from
Corollary \ref{surjective}, that
$\varphi_+$ is injective and that it is
actually a bijection onto $W_+$. This fact justifies the name \textit{biologically meaningful} points given to the points
of $W_+$. The
following result tells us that if $q\in
W_{\mathbf{\Delta}}$ has only positive
coordinates, then $q$ is non-singular
and $|\varphi^{-1}(q)|=4^{n-2}$. Among
all these pre-images, just the single
one in $\prod_{e\in E(T)}\D^{3}_+$
has biological meaning.
\end{rem}

\begin{cor} \label{singular_points}
A point $q=\varphi_0(p)\in W_{\D}$ is  singular if and only if
there is some $e\in E(T)$ such that $P^{(e)}\in \D^{3} \cap
\{P^e_C P^e_G P^e_T=0\}$. In particular, no point with biological
meaning is singular.
\end{cor}

\begin{proof}
It is well-known that the singular points on $W$ are the points
$\{q\in W \mid |\varphi^{-1}(q)|<4^{n-2}\}$, i.e. those points for
which at least one of their pre-images are invariant by the action
of some $g\in G'$. By Remark \ref{reflection}, we know that these
are precisely the points lying on an axis $P^e_x=0$ for some $e\in
E(T)$ and some $x\neq A$. As a consequence, the points in $W_+=\varphi(\prod_{e\in E(T)}\D_+^3)$ are non-singular. 
\end{proof}

\section{Local Complete Intersection}

Given a tree with $n$ leaves, the main
purpose of this section is to describe
a procedure to determine a local
complete intersection equal to the
variety $W=W_n$ in the open set $\D_+$.

\begin{notation}
 For every $n\in \mathbb{N}$, we will
 write $c(n)=4^{n-1}-6n+9$. Note that in virtue
 of Corollary \ref{dim_codim},
 $c(n)$ is the codimension of
 $W_n$.
\end{notation}

In Corollary \ref{singular_points} we
have seen that the points of $(W_n)_+$
are non-singular. Since any
regular local ring is a complete
intersection, the variety $W_n$ is a
\emph{local complete intersection} at
these points, i.e. the ideal $I_W$ can
be generated by $c(n)$ polynomials in a
neighborhood of these points or more
precisely, a minimal system of
generators for the localization of
$I_W$ in these points consists of
$c(n)=4^{n-1}-6n+9$ elements.

The following lemma provides a minimal
system of generators for this ideal in
case the tree $T$ has $n=3$ leaves.

\begin{prop}\label{gen-3leaves}
Let $T$ be a tree with 3 leaves and let
$W_3\subset \C^{16}$ be the model
associated to it. Then, the set of
quartics
\begin{eqnarray*}
h_1=q_{AAA}q_{ATT}q_{TCG}q_{TGC}-q_{ACC}q_{AGG}q_{TAT}q_{TTA}\\
h_2=q_{CCA}q_{CTG}q_{TAT}q_{TGC}-q_{CAC}q_{CGT}q_{TCG}q_{TTA}\\
h_3=q_{AGG}q_{ATT}q_{CAC}q_{CCA}-q_{AAA}q_{ACC}q_{CGT}q_{CTG}\\
h_4=q_{ACC}q_{ATT}q_{GAG}q_{GGA}-q_{AAA}q_{AGG}q_{GCT}q_{GTC}\\
h_5=q_{CAC}q_{CTG}q_{GCT}q_{GGA}-q_{CCA}q_{CGT}q_{GAG}q_{GTC}\\
h_6=q_{GGA}q_{GTC}q_{TAT}q_{TCG}-q_{GAG}q_{GCT}q_{TGC}q_{TTA}
\end{eqnarray*}
together with the equation $h=q_{AAA}-1$
is a local minimal system of generators
for the ideal of $W_3$ at the points of
$(W_3)_+$. Namely,
$\{h_1,h_2,h_3,h_4,h_5,h_6,q_{AAA}-1\}$
generate the ideal $I_{W_3}$ in the local ring
$\OO_{W,q}$, for any $q\in (W_3)_+$.
\end{prop}

\begin{rem}\label{nodepen}
It is worth pointing out that the minimal system of generators given in Proposition \ref{gen-3leaves} does not depend on the point $q \in (W_3)_+$. In the same way, the local complete intersection we will construct for an arbitrary tree of $n$ leaves will be the same for all points in $(W_n)_+$.
\end{rem}

\vspace{3mm}
\noindent \textsc{Proof of Proposition \ref{gen-3leaves}.}
Let $W'\subset \C^{16}$ be the variety
defined by the zero set of the ideal
$(h_1,h_2,h_3,h_4,h_5,h_6,h=q_{AAA}-1)$
and let $q\in (W_3)_+$. Let
$Jac_{W'}(q)$ be the jacobian matrix of
$W'$ at $q$:

\begin{tiny}
$\hspace*{6mm}  {q_{AAA}} \;
\hspace*{1mm} {q_{ACC}} \;
\hspace*{1mm} {q_{AGG}}
 \; \hspace*{1mm}  {q_{ATT}} \; \hspace*{1mm} {q_{CAC}} \; \hspace*{1mm} {q_{CCA}}
 \; \hspace*{1mm}  {q_{CGT}} \; \hspace*{1mm} {q_{CTG}} \; \hspace*{1mm} {q_{GAG}}
 \; \hspace*{1mm}  {q_{GCT}} \; \hspace*{1mm} {q_{GGA}} \; \hspace*{1mm} {q_{GTC}}
 \; \hspace*{1mm}  {q_{TAT}} \; \hspace*{1mm} {q_{TCG}} \; \hspace*{1mm} {q_{TGC}}
 \; \hspace*{1mm}  {q_{TTA}}$
\end{tiny}~\begin{small}\begin{eqnarray*}
\begin{array}{c}
{h_1} \\
{h_2} \\
{h_3} \\
{h_4} \\
{h_5} \\
{h_6} \\
{h}
\end{array}
\left (
\begin{array}
{cccccccccccccccc}
* \quad  & * \quad  & * \quad & * \quad & 0 \quad & 0 \quad & 0 \quad & 0 \quad & 0 \quad & 0 \quad & 0 \quad & 0 \quad & * \quad & * \quad & * \quad & * \quad \\
0 \quad  & 0 \quad  & 0 \quad & 0 \quad & * \quad & * \quad & * \quad & * \quad & 0 \quad & 0 \quad & 0 \quad & 0 \quad & * \quad & * \quad & * \quad & * \quad \\
* \quad  & * \quad  & * \quad & * \quad & * \quad & * \quad & * \quad & * \quad & 0 \quad & 0 \quad & 0 \quad & 0 \quad & 0 \quad & 0 \quad & 0 \quad & 0 \quad \\
* \quad  & * \quad  & * \quad & * \quad & 0 \quad & 0 \quad & 0 \quad & 0 \quad & * \quad & * \quad & * \quad & * \quad & 0 \quad & 0 \quad & 0 \quad & 0 \quad \\
0 \quad  & 0 \quad  & 0 \quad & 0 \quad & * \quad & * \quad & * \quad & * \quad & * \quad & * \quad & * \quad & * \quad & 0 \quad & 0 \quad & 0 \quad & 0 \quad \\
0 \quad  & 0 \quad  & 0 \quad & 0 \quad & 0 \quad & 0 \quad & 0 \quad & 0 \quad & * \quad & * \quad & * \quad & * \quad & * \quad & * \quad & * \quad & * \quad \\
1 \quad  & 0 \quad  & 0 \quad & 0 \quad
& 0 \quad & 0 \quad & 0 \quad & 0 \quad
& 0 \quad & 0 \quad & 0 \quad & 0 \quad
& 0 \quad & 0 \quad & 0 \quad & 0 \quad
\end{array}
\right )
\end{eqnarray*}
\end{small}
with entries
\begin{eqnarray*}
\begin{array}{ll}
(h_1,q_{AAA})=q_{ATT}q_{TCG}q_{TGC} &
(h_1,q_{ACC})=-q_{AGG}q_{TAT}q_{TTA}\\
(h_1,q_{AGG})=-q_{ACC}q_{TAT}q_{TTA} &
(h_1,q_{ATT})=q_{AAA}q_{TCG}q_{TGC}\\
(h_1,q_{TAT})=-q_{ACC}q_{AGG}q_{TTA} &
(h_1,q_{TCG})=q_{AAA}q_{ATT}q_{TGC}\\
(h_1,q_{TGC})=q_{AAA}q_{ATT}q_{TCG} &
(h_1,q_{TTA})=-q_{ACC}q_{AGG}q_{TAT}
\end{array}
\end{eqnarray*}
\begin{eqnarray*}
\begin{array}{ll}
(h_2,q_{CAC})=-q_{CGT}q_{TCG}q_{TTA}&
(h_2,q_{CCA})=q_{CTG}q_{TAT}q_{TGC}\\
(h_2,q_{CGT})=-q_{CAC}q_{TCG}q_{TTA}&
(h_2,q_{CTG})=q_{CCA}q_{TAT}q_{TGC}\\
(h_2,q_{TAT})=q_{CCA}q_{CTG}q_{TGC}&
(h_2,q_{TCG})=-q_{CAC}q_{CGT}q_{TTA}\\
(h_2,q_{TGC})=q_{CCA}q_{CTG}q_{TAT}&
(h_2,q_{TTA})=-q_{CAC}q_{CGT}q_{TCG}
\end{array}
\end{eqnarray*}
\begin{eqnarray*}
\begin{array}{ll}
(h_3,q_{AAA})=-q_{ACC}q_{CGT}q_{CTG}&
(h_3,q_{ACC})=-q_{AAA}q_{CGT}q_{CTG}\\
(h_3,q_{AGG})=q_{ATT}q_{CAC}q_{CCA}&
(h_3,q_{ATT})=q_{AGG}q_{CAC}q_{CCA}\\
(h_3,q_{CAC})=q_{AGG}q_{ATT}q_{CCA}&
(h_3,q_{CCA})=q_{AGG}q_{ATT}q_{CCA}\\
(h_3,q_{CGT})=-q_{AAA}q_{ACC}q_{CTG}&
(h_3,q_{CTG})=-q_{AAA}q_{ACC}q_{CGT}
\end{array}
\end{eqnarray*}
\begin{eqnarray*}
\begin{array}{ll}
(h_4,q_{AAA})=-q_{AGG}q_{GCT}q_{GTC}&
(h_4,q_{ACC})=q_{ATT}q_{GAG}q_{GGA}\\
(h_4,q_{AGG})=-q_{AAA}q_{GCT}q_{GTC}&
(h_4,q_{ATT})=q_{ACC}q_{GAG}q_{GGA}\\
(h_4,q_{GAG})=q_{ACC}q_{ATT}q_{GGA}&
(h_4,q_{GCT})=-q_{AAA}q_{AGG}q_{GTC}\\
(h_4,q_{GGA})=q_{ACC}q_{ATT}q_{GAG}&
(h_4,q_{GTC})=-q_{AAA}q_{AGG}q_{GCT}
\end{array}
\end{eqnarray*}
\begin{eqnarray*}
\begin{array}{ll}
(h_5,q_{CAC})=q_{CTG}q_{GCT}q_{GGA}&
(h_5,q_{CCA})=-q_{CGT}q_{GAG}q_{GTC}\\
(h_5,q_{CGT})=-q_{CCA}q_{GAG}q_{GTC}&
(h_5,q_{CTG})=q_{CAC}q_{GCT}q_{GGA}\\
(h_5,q_{GAG})=-q_{CCA}q_{CGT}q_{GTC}&
(h_5,q_{GCT})=q_{CAC}q_{CTG}q_{GGA}\\
(h_5,q_{GGA})=q_{CAC}q_{CTG}q_{GCT}&
(h_5,q_{GTC})=-q_{CCA}q_{CGT}q_{GAG}
\end{array}
\end{eqnarray*}
\begin{eqnarray*}
\begin{array}{ll}
(h_6,q_{GAG})=-q_{GCT}q_{TGC}q_{TTA}&
(h_6,q_{GCT})=-q_{GAG}q_{TGC}q_{TTA}\\
(h_6,q_{GGA})=q_{GTC}q_{TAT}q_{TCG}&
(h_6,q_{GTC})=q_{GGA}q_{TAT}q_{TCG}\\
(h_6,q_{TAT})=q_{GGA}q_{GTC}q_{TCG}&
(h_6,q_{TCG})=q_{GGA}q_{GTC}q_{TAT}\\
(h_6,q_{TGC})=-q_{GAG}q_{GCT}q_{TTA}&
(h_6,q_{TTA})=-q_{GAG}q_{GCT}q_{TGC}
\end{array}
\end{eqnarray*}

In general, one has that
$rk(Jac_{q}(W'))\leq codim_{q}(W') \leq
7$. It can be seen by direct
computation that the $6\times 6$-matrix
obtained from $Jac_{W'}(q)$ by removing
the last row and keeping the columns
indexed by $q_{ACC}$, $q_{ATT}$,
$q_{CAC}$, $q_{CTG}$, $q_{TAT}$,
$q_{TCG}$ equals
\begin{tiny}
\begin{eqnarray*}
2q_{GAG}q_{GCT}q_{GGA}^3q_{GTC} & ( & q_{AAA}q_{AGG}q_{ATT}^2q_{CAC}q_{CCA}^2q_{CTG}q_{TAT}q_{TCG}q_{TGC}^2+ \\
q_{AAA}^2q_{ACC}q_{ATT}q_{CCA}q_{CGT}q_{CTG}^2q_{TAT}q_{TCG}q_{TGC}^2 & + &
q_{ACC}q_{AGG}^2q_{ATT}q_{CAC}q_{CCA}^2q_{CTG}q_{TAT}^2q_{TGC}q_{TTA}+\\
q_{AAA}q_{ACC}^2q_{AGG}q_{CCA}q_{CGT}q_{CTG}^2q_{TAT}^2q_{TGC}q_{TTA}& + &
q_{AAA}q_{AGG}q_{ATT}^2q_{CAC}^2q_{CCA}q_{CGT}q_{TCG}^2q_{TGC}q_{TTA}+\\
q_{AAA}^2q_{ACC}q_{ATT}q_{CAC}q_{CGT}^2q_{CTG}q_{TCG}^2q_{TGC}q_{TTA}& + &
q_{ACC}q_{AGG}^2q_{ATT}q_{CAC}^2q_{CCA}q_{CGT}q_{TAT}q_{TCG}q_{TTA}^2+\\
q_{AAA}q_{ACC}^2q_{AGG}q_{CAC}q_{CGT}^2q_{CTG}q_{TAT}q_{TCG}q_{TTA}^2 & ) &
\end{eqnarray*}
\end{tiny}which is clearly positive in $\D_+$.
Therefore,  $rk(Jac_{q}(W'))=7$ and so,
$q$ is a non-singular point of $W'$ and
$W'$ is a local complete intersection
at $q$. Therefore, $W' \subset
\CC^{16}$  is a subvariety of dimension
9 containing $W_3$, which has also
dimension 9 and is non-singular at $q$.
It follows that $W'$ and $W_3$ coincide
in a neighborhood of $q$ and we are
done. \qed

\begin{rem}\label{induction}
For future reference, it is worth
noting that the matrix $J'$ obtained
from $Jac_q(W')$ by removing the
columns
$q_{AAA},q_{CCA},q_{GGA},q_{TTA}$ and
the last row has maximal rank equal to
6.
\end{rem}

Next, we want to describe a procedure
to give a minimal system of generators
for the ideal of $W_n$ around any point
$q\in (W_n)_+$. Some of these generators are
determined recursively from subtrees of
$T$, while the remaining are easily
inferred from some matrices to be
defined later.

\vspace{3mm}

First we describe how these
generators are to be constructed by induction on the number of leaves. Then,
we will prove that the whole set of
these polynomials generate a complete
intersection which equals the variety
$W_n$ in a neighborhood of any $q\in
(W_n)_+$. The generators of this local complete intersection ideal will not depend on the point $q$, as we pointed out in Remark \ref{nodepen}.

\subsection*{Generators of degree 4} As above, write $R=\mathbb{C}[\Q_n]$
for the ring of polynomials in the
unknowns $\Q_n$. Following the idea of
Chang \cite{Chang96}, write
$v_1,\ldots,v_n$ for the leaves of $T$.
By reordering the leaves, we may assume
that $v_{n-1}$ and $v_n$ form a
\textit{cherry}, i.e. are joined to a
node $m$. Take the tree $T'$ with
leaves $L(T')=L(T)\cup \{m\}-
\{v_{n-1},v_n\}$, interior nodes
$N(T')=N(T)-\{m\}$ and edges
$E(T')=E(T)-\{[m,v_{n-1}],[m,v_n]\}$,
where $[m,n]$ is the edge containing
the nodes $m$ and $n$ (see figure
\ref{demoinvariantsCI}). In virtue of
Corollaries \ref{surjective} and
\ref{dim_codim}, the variety $W_{n-1}$
associated to $T'$ is the image of the
polynomial map in (\ref{paramfourier})
\begin{eqnarray*}
\varphi_{n-1}: (\C^3 )^{2n-5} \lra
\C^{4^{n-2}}
\end{eqnarray*}
and has dimension $3(2n-5)$.
%

Assume that we have constructed a local complete intersection $\{g_1,\ldots,g_{c(n-1)}\}$ at the
points of $(W_{n-1})_+$
(equivalently,
$\{g_1,\ldots,g_{c(n-1)}\}$ generate
the localization of the ideal $I_{T'}$
at the points of $(W_{n-1})_+$).
The map $j_{n-1}:\Q_{n-1}\rightarrow
\Q_n$ defined by $q_{x_1\ldots x_{n-1}}
\mapsto q_{x_1\ldots x_{n-1}A}$ induces
a ring homomorphism
\begin{eqnarray*}
\psi_{n-1}:\C[\Q_{n-1}] & \rightarrow & R.
\end{eqnarray*}
Write
\begin{eqnarray*}\label{1stgroup}
J(n-1)=\{f^{(n-1)}_1,\ldots,f^{(n-1)}_{c(n-1)}\}\subset
R
\end{eqnarray*}
for the set of polynomials being the
image by $\psi_{n-1}$ of the generators
$\{g_i\}$.

Analogously, let $T''$ be the tree with
3 leaves determined by the vertices
$v_1,v_{n-1}$ and $v_n$. The
variety $W_3\subset \C^{16}$ is the image of
\begin{eqnarray*} \varphi_3: (\C^4
)^{3} \lra \C^{16}
\end{eqnarray*}
and has dimension 9. A complete
system of generators
$\{h_1,\ldots,h_6\}$ of the ideal
$I_{T''}\subset \C[\Q_3]$ is given by
Lemma \ref{gen-3leaves}. As above, the
map $j_3:\Q_3\rightarrow \Q_n$ defined
by $q_{xyz} \mapsto q_{xA\ldots Ayz}$
induces a ring homomorphism
\begin{eqnarray*}
\psi_3:\C[\Q_3] & \rightarrow & R.
\end{eqnarray*}
Write
\begin{eqnarray*}\label{2ndgroup}
J(3)=\{f^{(3)}_1,\ldots,f^{(3)}_6\}
\subset R
\end{eqnarray*}
for the set of polynomials being the
image by $\psi_{n-1}$ of $\{h_1,\ldots,h_6\}$.
The polynomials in $J(3)$ and $J(n-1)$ are quartics, but we still need to construct an extra set of polynomials of degree 2.

\subsection*{Generators of degree 2} Now, for each letter $z\in \Sigma$, write
$M(z)$ for the $4\times 4^{n-3}$-matrix
with rows indexed by the couples
$\{xy \mid x+y=z\}$, columns indexed by
$\{x_1\dots x_{n-2} \mid \sum_{i=1}^{n-2}
x_i=z\}$ and whose $(xy,x_1\dots
x_{n-2})$-entry is precisely $q_{x_1\dots
x_{n-2}xy}$:
\begin{eqnarray*}
& & \qquad \qquad \qquad \quad \; \scriptstyle{x_1\dots x_{n-2}} \\
M(z)& = & \scriptstyle{xy} \quad \left (
\begin{array} {ccccc} \ddots \quad &
& \vdots & & \quad \ddots \\
\cdots & & q_{x_1\dots x_{n-2}xy} & & \cdots \\
\ddots & & \vdots & & \ddots
\end{array} \right )
\end{eqnarray*}
For each of these matrices, take
the set of the $3(4^{n-3}-1)$ $2\times
2$-minors containing $q_{zA\ldots
AzA}$: we obtain
polynomials $F^{(z)}_i\in R$ of the
form
\begin{eqnarray}\label{formK}
q_{x_1\ldots x_{n-2}xy}q_{zA\ldots
AzA}-q_{x_1\ldots x_{n-2}zA}q_{zA\ldots
Axy}=0
\end{eqnarray}
for $i=1,\ldots, 3(4^{n-3}-1)$. We get a
total of $12(4^{n-3}-1)$ polynomials.
 For each letter $z\in
\Sigma$, write $K(z)=\{F^{(z)}_i\}$ for
this set of polynomials and
\begin{eqnarray*}\label{3rdgroup}
K=\bigcup_{z\in \Sigma} K(z).
\end{eqnarray*}

\begin{thm}\label{invariantsCI}
At each point $q \in (W_n)_+$, the ideal generated by the set
\[J(3)\cup J(n-1) \cup K\]
together with the equation $q_{A\ldots
A}=1$ is a local complete intersection that defines the Kimura variety $W_n$ in a neighborhood of $q$.
\end{thm}

\begin{proof}
First of all, direct computation shows
that the number of polynomials being
considered equals the  codimension of
the variety $W$, i.e.
\begin{eqnarray*}
|J(3)|+|J(n-1)| + |K|+1  & = & \\
6+(4^{n-2}-6(n-1)+8)+12(4^{n-3}-1)+1
& = & 4^{n-1}-6n
+9.
\end{eqnarray*}
By \cite{Sturmfels2005} we know that
the  ideal $\wp=(q_{A\ldots
A}-1,J(3),J(n-1),K)$ is contained in $I_{W_n}$ and
so, $W_n$ is contained in the variety
$W'$ defined by $\wp$. We claim that $q$ is non-singular
in $W'$. From this, we deduce that $W'$
is a local complete intersection at
$q$, and as we did in the proof of
Lemma \ref{gen-3leaves}, we conclude
that it is equal to $W_n$ in a
neighborhood of $q$.

Now we prove that $q$ is a smooth point
of $W'$. Write $Q_0=\{q_{AA\ldots
AAA},q_{CA\ldots ACA},\\q_{GA\ldots
AGA},q_{TA\ldots ATA}\}$, and notice
that \[Q_0=j_3(\Q_3)\cap
j_{n-1}(\Q_{n-1}).\] By reordering the
rows  and columns if necessary, we may
assume that the jacobian matrix of $W'$
at $q$ has the form
\begin{eqnarray*}
Jac_{q}(W')=\left (
\begin{array}
{ccc}
  B & 0 & 0 \\
  0 & J' & 0 \\
  * & * & D
\end{array}
\right )
\end{eqnarray*}
where the $c(n-1)\times 4^{n-2}$-matrix
$B$ equals the jacobian matrix
$Jac_{q_{n-1}}(W_{n-1})$, $J'$ is the
$6\times 12$-matrix of Lemma
\ref{induction}. In this way, the
columns of the submatrix $B$ are
indexed by the unknowns in $\Q_{n-1}$
while the rows are indexed by the
equations $\{q_{A\ldots A}-1\}\cup
J(n-1)$. Similarly, the columns of $J'$
are indexed by $\Q_3\setminus
\{q_{AAA},q_{CCA},q_{GGA},q_{TTA}\}$
while the rows are indexed by $J(3)$.
The columns of the matrix $D$ are
indexed by the remaining unknowns while
its rows are indexed by the equations
$\{F_i^{(z)}\}_{z\in \Sigma}$. Each
of these equations has the form
(\ref{formK}) and so, its partial
derivative relative to the unknown
$q_{x_1\ldots x_{n-2}xy}$ is equal to
$q_{zA\ldots zA}$. Therefore, by
reordering rows and columns if
necessary we may assume that the matrix
$D$ is a diagonal matrix (and all its entries are strictly positive because $q \in (W_n)_+$).

In virtue of \ref{induction}, we  know
that $rank(J')=6$ and, by induction
hypothesis, $B$ has maximal rank equal
to $c(n-1)=4^{n-2}-6n+15$. It follows
that the matrix $Jac_{q}(W')$ has maximal
rank equal to
\begin{eqnarray*}
rank(Jac_{q}(W'))=4^{n-2}-6n+15 + 6 +12(4^{n-3}-1) = 4^{n-1}-6n+9,
\end{eqnarray*}
and we are done.
\begin{figure}
\begin{center}
\psfrag{a}{$v_{n}$}\psfrag{b}{$v_{n-1}$}\psfrag{c}{$v_1$}\psfrag{m}{$m$}
\psfrag{T}{$T'$}

\includegraphics[scale=0.4]{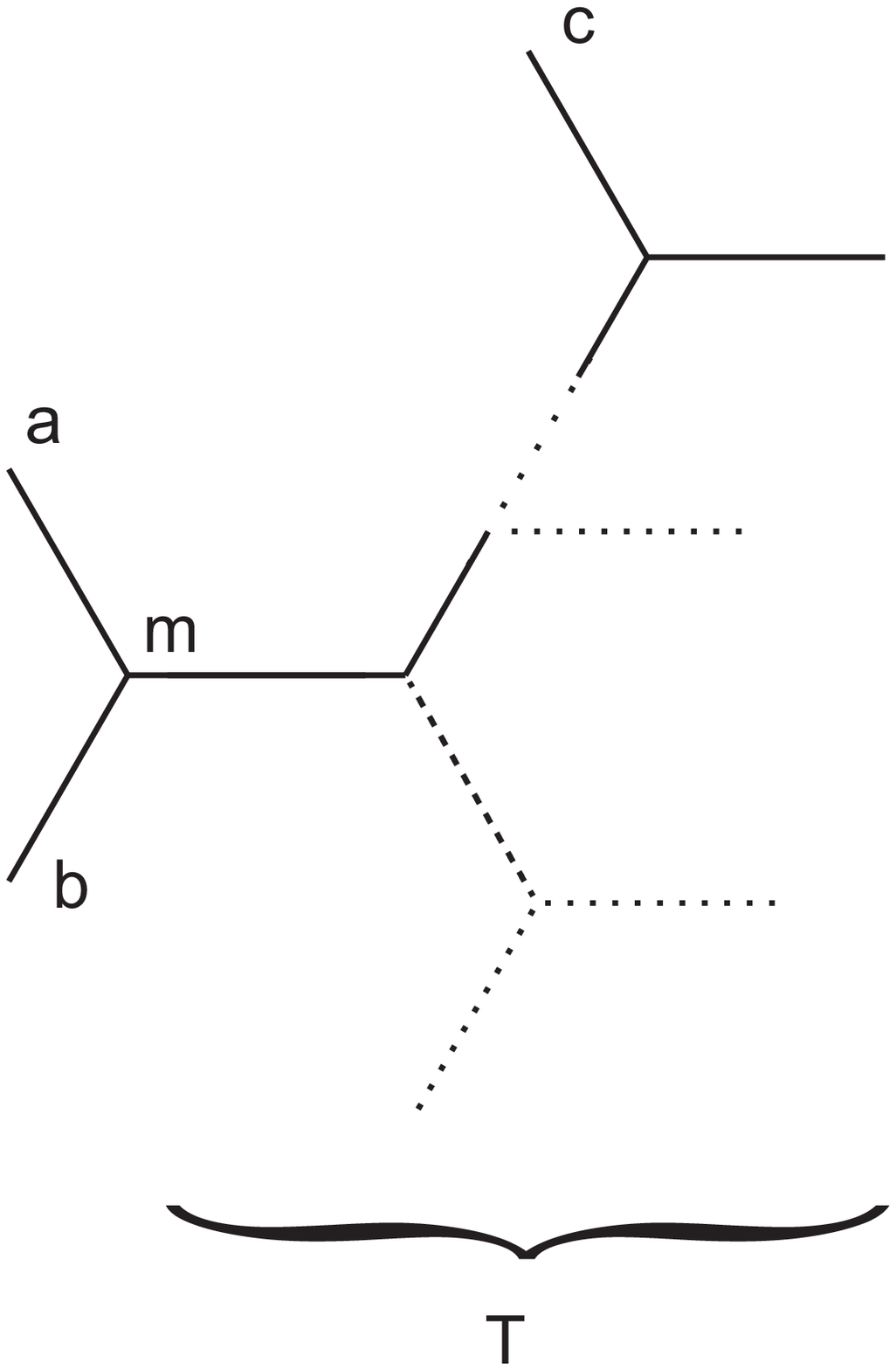}
\end{center}
\caption{\label{demoinvariantsCI}}
\end{figure}
\end{proof}

\begin{rem}
The set of quadrics $K$ contains the information of invariants coming from the splits of the tree (see \cite{Sturmfels2005} and \cite{Eriksson05}). Although in theory a tree can be reconstructed from its splits (see \cite{Eriksson05} Theorem 19.14), the variety defined by $K$ is much bigger than $W$ because it has codimension $12(4^{n-3}-1)$.
\end{rem}
\begin{rem}
In \cite{CFS}, we studied  a
phylogenetic reconstruction method
(already introduced in \cite{CGS})
which was based on a set of generators
of the ideal associated to the Kimura
model. The simulation studies performed
there showed that it is actually a very
competitive and highly efficient
method. In the case of 4-leaved trees
and for the Kimura 3-parameter model, a
minimal system of generators for the
corresponding ideal consists of 8002
polynomials of degrees 2, 3 and 4.
Because of the results of this paper,
it is enough to deal with the 48
invariants listed in the following
example (or in general, the codimension
of the variety $W_n$). This leads to a
substantial improvement in the
efficiency and effectiveness of the
method. Simulations studies on this
variant of the method can be seen on
the webpage
\begin{center}
http://wwww.ma1.upc.edu/$\scriptstyle \sim$jfernandez/ci.html
\end{center}
and the reader should contrast them to \cite{CFS}. 
Moreover, the fact that we provide the \textit{smallest} set of local generators  in Theorem \ref{invariantsCI}, gives some hope for the generalization of phylogenetic reconstruction methods based on algebraic geometry to trees with a large number of leaves.

\end{rem}

\begin{example}\label{4leavestree}
Let $T$ be the unrooted 4-leaved tree of figure \ref{clawtree}(b). The above procedure gives rise to the following 48 invariants:
36 quadrics
\begin{small}
\begin{eqnarray*}
q_{CCCC}q_{AAAA}-q_{CCAA}q_{AACC}, & q_{GGCC}q_{AAAA}-q_{GGAA}q_{AACC}, & q_{TTCC}q_{AAAA}-q_{TTAA}q_{AACC}, \\
q_{CCGG}q_{AAAA}-q_{CCAA}q_{AAGG}, & q_{GGGG}q_{AAAA}-q_{GGAA}q_{AAGG}, & q_{TTGG}q_{AAAA}-q_{TTAA}q_{AAGG}, \\
q_{CCTT}q_{AAAA}-q_{CCAA}q_{AATT}, & q_{GGTT}q_{AAAA}-q_{GGAA}q_{AATT}, & q_{TTTT}q_{AAAA}-q_{TTAA}q_{AATT}, \\
q_{ACAC}q_{CACA}-q_{ACCA}q_{CAAC}, & q_{GTAC}q_{CACA}-q_{GTCA}q_{CAAC}, & q_{TGAC}q_{CACA}-q_{TGCA}q_{CAAC}, \\
q_{ACGT}q_{CACA}-q_{ACCA}q_{CAGT}, & q_{GTGT}q_{CACA}-q_{GTCA}q_{CAGT}, & q_{TGGT}q_{CACA}-q_{TGCA}q_{CAGT}, \\
q_{ACTG}q_{CACA}-q_{ACCA}q_{CATG}, & q_{GTTG}q_{CACA}-q_{GTCA}q_{CATG}, & q_{TGTG}q_{CACA}-q_{TGCA}q_{CATG}, \\
q_{AGAG}q_{GAGA}-q_{AGGA}q_{GAAG}, & q_{CTAG}q_{GAGA}-q_{CTGA}q_{GAAG}, & q_{TCAG}q_{GAGA}-q_{TCGA}q_{GAAG}, \\
q_{AGCT}q_{GAGA}-q_{AGGA}q_{GACT}, & q_{CTCT}q_{GAGA}-q_{CTGA}q_{GACT}, & q_{TCCT}q_{GAGA}-q_{TCGA}q_{GACT}, \\
q_{AGTC}q_{GAGA}-q_{AGGA}q_{GATC}, & q_{CTTC}q_{GAGA}-q_{CTGA}q_{GATC}, & q_{TCTC}q_{GAGA}-q_{TCGA}q_{GATC}, \\
q_{ATAT}q_{TATA}-q_{ATTA}q_{TAAT}, & q_{CGAT}q_{TATA}-q_{CGTA}q_{TAAT}, & q_{GCAT}q_{TATA}-q_{GCTA}q_{TAAT}, \\
q_{ATCG}q_{TATA}-q_{ATTA}q_{TACG}, & q_{CGCG}q_{TATA}-q_{CGTA}q_{TACG}, & q_{GCCG}q_{TATA}-q_{GCTA}q_{TACG}, \\
q_{ATGC}q_{TATA}-q_{ATTA}q_{TAGC}, & q_{CGGC}q_{TATA}-q_{CGTA}q_{TAGC}, & q_{GCGC}q_{TATA}-q_{GCTA}q_{TAGC},
\end{eqnarray*}
\end{small}
and 12 quartics
\begin{small}
\begin{eqnarray*}
& q_{AAAA}q_{ATTA}q_{TCGA}q_{TGCA}-q_{ACCA}q_{AGGA}q_{TATA}q_{TTAA}, \\
& q_{CCAA}q_{CTGA}q_{TATA}q_{TGCA}-q_{CACA}q_{CGTA}q_{TCGA}q_{TTAA}, \\
& q_{AGGA}q_{ATTA}q_{CACA}q_{CCAA}-q_{AAAA}q_{ACCA}q_{CGTA}q_{CTGA}, \\
& q_{ACCA}q_{ATTA}q_{GAGA}q_{GGAA}-q_{AAAA}q_{AGGA}q_{GCTA}q_{GTCA}, \\
\end{eqnarray*}
\begin{eqnarray*}
& q_{CACA}q_{CTGA}q_{GCTA}q_{GGAA}-q_{CCAA}q_{CGTA}q_{GAGA}q_{GTCA}, \\
& q_{GGAA}q_{GTCA}q_{TATA}q_{TCGA}-q_{GAGA}q_{GCTA}q_{TGCA}q_{TTAA}, \\
& q_{AAAA}q_{AATT}q_{TACG}q_{TAGC}-q_{AACC}q_{AAGG}q_{TAAT}q_{TATA}, \\
& q_{CACA}q_{CATG}q_{TAAT}q_{TAGC}-q_{CAAC}q_{CAGT}q_{TACG}q_{TATA}, \\
& q_{AAGG}q_{AATT}q_{CAAC}q_{CACA}-q_{AAAA}q_{AACC}q_{CAGT}q_{CATG}, \\
& q_{AACC}q_{AATT}q_{GAAG}q_{GAGA}-q_{AAAA}q_{AAGG}q_{GACT}q_{GATC}, \\
& q_{CAAC}q_{CATG}q_{GACT}q_{GAGA}-q_{CACA}q_{CAGT}q_{GAAG}q_{GATC}, \\
& q_{GAGA}q_{GATC}q_{TAAT}q_{TACG}-q_{GAAG}q_{GACT}q_{TAGC}q_{TATA}.
\end{eqnarray*}
\end{small}
\end{example}


\begin{thebibliography}{SSEW93}

\bibitem[AR03]{Allman2003}
ES~Allman and JA~Rhodes.
\newblock Phylogenetic invariants for the general {M}arkov model of sequence
  mutation.
\newblock {\em Math. Biosci.}, 186(2):113--144, 2003.

\bibitem[AR04]{ARbook}
ES~Allman and JA~Rhodes.
\newblock {\em Mathematical models in biology, an introduction}.
\newblock Cambridge University Press, January 2004.
\newblock ISBN 0-521-52586-1).

\bibitem[AR07]{Allman2004b}
ES~Allman and JA~Rhodes.
\newblock Phylogenetic ideals and varieties for the general {M}arkov model.
\newblock to appear in Advances in Applied Mathematics, 2007.

\bibitem[CFS07]{CFS}
M~Casanellas and J~Fernandez-Sanchez.
\newblock Performance of a new invariants method on homogeneous and
  nonhomogeneous quartet trees.
\newblock {\em Molecular Biology and Evolution}, 24(1):288--293, January 2007.

\bibitem[CGS05]{CGS}
M~Casanellas, LD~Garcia, and S~Sullivant.
\newblock Catalog of small trees.
\newblock In L.~Pachter and B.~Sturmfels, editors, {\em Algebraic Statistics
  for computational biology}, chapter~15. Cambridge University Press, 2005.

\bibitem[Cha96]{Chang96}
J~T Chang.
\newblock Full reconstruction of {M}arkov models on evolutionary trees:
  identifiability and consistency.
\newblock {\em Math. Biosci.}, 137(1):51--73, 1996.

\bibitem[CS05]{CS}
M~Casanellas and S~Sullivant.
\newblock The strand symmetric model.
\newblock In L.~Pachter and B.~Sturmfels, editors, {\em Algebraic Statistics
  for computational biology}, chapter~16. Cambridge University Press, 2005.

\bibitem[Dol03]{Dolgachev}
I~Dolgachev.
\newblock {\em Lectures on invariant theory}, volume 296 of {\em London
  Mathematical Society Lecture Note Series}.
\newblock Cambridge University Press, Cambridge, 2003.

\bibitem[Eri05]{Eriksson05}
N~Eriksson.
\newblock Tree construction using singular value decomposition.
\newblock In L.~Pachter and B.~Sturmfels, editors, {\em Algebraic Statistics
  for computational biology}, chapter~19, pages 347--358. Cambridge University
  Press, 2005.

\bibitem[ERSS05]{PhyloAG}
N~Eriksson, K~Ranestad, B~Sturmfels, and S~Sullivant.
\newblock Phylogenetic algebraic geometry.
\newblock In C~Ciliberto, A~Geramita, B~Harbourne, R-M Roig, and K~Ranestad,
  editors, {\em Projective Varieties with Unexpected Properties}, pages
  237--255. De Gruyter, 2005.

\bibitem[ES93]{Evans1993}
S~Evans and T~Speed.
\newblock Invariants of some probability models used in phylogenetic inference.
\newblock {\em The Annals of Statistics}, 21:355--377, 1993.

\bibitem[Hag00]{Hagedorn2000}
Thomas~R. Hagedorn.
\newblock Determining the number and structure of phylogenetic invariants.
\newblock {\em Adv. in Appl. Math.}, 24:1--21, 2000.

\bibitem[JC69]{JC69}
T.~H. Jukes and C.~R. Cantor.
\newblock Evolution of protein molecules.
\newblock {\em In Mammalian Protein Metabolism}, pages 21--132, 1969.

\bibitem[Kim80]{Kimura1980}
M~Kimura.
\newblock A simple method for estimating evolutionary rates of base
  substitution through comparative studies of nucleotide sequences.
\newblock {\em J. Mol. Evol.}, 16:111--120, 1980.

\bibitem[Kim81]{Kimura1981}
M~Kimura.
\newblock Estimation of evolutionary sequences between homologous nucleotide
  sequences.
\newblock {\em Proc. Nat. Acad. Sci. , USA}, 78:454--458, 1981.

\bibitem[MS05]{MillerSturmfels}
E~Miller and B~Sturmfels.
\newblock {\em Combinatorial commutative algebra}, volume 227 of {\em Graduate
  Texts in Mathematics}.
\newblock Springer-Verlag, New York, 2005.

\bibitem[PS04]{Pachter2004}
L~Pachter and B~Sturmfels.
\newblock Tropical geometry of statistical models.
\newblock {\em Proceedings of the National Academy of Sciences},
  101:16132--16137, 2004.

\bibitem[PS05]{ASCB2005}
L~Pachter and B~Sturmfels, editors.
\newblock {\em Algebraic Statistics for computational biology}.
\newblock Cambride University Press, November 2005.
\newblock ISBN 0-521-85700-7.

\bibitem[SHP98]{Steel98}
MA~Steel, MD~Hendy, and D~Penny.
\newblock Reconstructing phylogenies from nucleotide pattern probabilities: a
  survey and some new results.
\newblock {\em Discrete Appl. Math.}, 88(1-3):367--396, 1998.

\bibitem[SS05]{Sturmfels2005}
B~Sturmfels and S~Sullivant.
\newblock Toric ideals of phylogenetic invariants.
\newblock {\em J. Comput. Biol.}, 12:204--228, 2005.

\bibitem[SSEW93]{Steel1993}
MA~Steel, LA~Sz{\'e}kely, PL~Erd{\H{o}}s, and P~Waddell.
\newblock A complete family of phylogenetic invariants for any number of taxa
  under kimura's 3st model.
\newblock {\em NZ J. Bot.}, 31:289--296, 1993.

\bibitem[WB06]{Wisniewski}
JA~Wisniewski and W~Buczynska.
\newblock On phylogenetic trees - a geometer's view.
\newblock http://arxiv.org/abs/math.AG/0601357, 2006.

\end{thebibliography}
\end{document}